\documentclass[a4paper,twoside,12pt]{article}
\usepackage{amsmath,amssymb,amsthm}
\usepackage{verbatim}
\usepackage{bbm}
\usepackage{stmaryrd}
\usepackage{graphicx,color}
\usepackage{float}
\usepackage{caption}
\usepackage{pifont}
\usepackage{wasysym}
\usepackage{framed}
\usepackage{indentfirst}

\usepackage{placeins}

\usepackage{enumitem}
\setlist{itemsep=0pt,topsep=2pt,partopsep=2pt}

\usepackage{xspace}

\usepackage{titlesec}

\titleformat*{\section}{\large\bfseries}
\titleformat*{\subsection}{\normalsize\bfseries}

\usepackage[table,dvipsnames]{xcolor}
\newcommand{\mostimportant}[1]{{\footnotesize #1\par}}

\DeclareMathOperator{\SO}{SO}
\DeclareMathOperator{\OO}{O}
\DeclareMathOperator{\tr}{tr}

\DeclareMathOperator{\ANGLE}{angle}

\newcommand{\dd}{\mathrm{d}}
\newcommand{\p}{\partial}

\renewcommand{\epsilon}{USE varepsilon INSTEAD}

\newcommand{\C}{\mathbbm{C}}
\newcommand{\Q}{\mathbbm{Q}}
\newcommand{\N}{\mathbbm{N}}
\newcommand{\R}{\mathbbm{R}}
\newcommand{\Z}{\mathbbm{Z}}

\newcommand{\step}{\vskip 3mm \noindent }

\newcommand{\KA}{\text{KA}}

\newcommand{\NN}[1]{\|#1\|_{\text{N}}}
\newcommand{\NDN}[1]{\|#1\|_{\text{DN}}}
\newcommand{\NNOP}[1]{\|#1\|_{\text{N} \leftarrow \text{N}}}

\newcommand{\SN}{\text{N}}
\newcommand{\SDN}{\text{DN}}
\newcommand{\SBN}{\text{BN}}

\newcommand{\SSS}{\text{\bf S}}

\newtheoremstyle{sproperty}%
 {20pt}% space above
 {3pt}% space below
 {}% body font
 {\parindent}% indent amount
 {}% head font
 {}% punctuation after theorem
 { }% space after theorem head
 {{\bf #1#2 (#3).}}% theorem head spec

\newcommand{\speculationtag}{P}
\theoremstyle{sproperty}
\newtheorem{sproperty}{\speculationtag}
\newcommand{\sref}[1]{\speculationtag\ref{#1}}

\newtheoremstyle{sstep}%
 {20pt}% space above
 {3pt}% space below
 {}% body font
 {\parindent}% indent amount
 {}% head font
 {}% punctuation after theorem
 { }% space after theorem head
 {{\bf #1#2.}}% theorem head spec

\newcommand{\steptag}{S}
\theoremstyle{sstep}
\newtheorem{sstep}{\steptag}
\newcommand{\stepref}[1]{\steptag\ref{#1}}

\begin{document}

%%%%%%%%%%%%%%%%%%%%%%%%%%%%%%%%%%%%%%%%%%%%%%%%%%%%%%%%%%%%%%%%%%%%%%%%%%%
\noindent{\bf\Large 3D incompressible Euler: Geometric formalism}\\
\noindent{\bf\Large\rule{0pt}{20pt}and a hypothetical self-similar flow}\\
\noindent\rule{0pt}{35pt}{\bf Michael Reiterer}
%%%%%%%%%%%%%%%%%%%%%%%%%%%%%%%%%%%%%%%%%%%%%%%%%%%%%%%%%%%%%%%%%%%%%%%%%%%
\vskip 13mm
\noindent {\bf Abstract:}
We give a geometric formulation of
3D incompressible Euler
that contains the Eulerian and Lagrangian gauges as special cases.
In the Lagrangian gauge, incompressible Euler is a
real analytic ODE in Banach space;
a short proof of this known result is given.
We then describe (in a self-contained section)
some basic properties of
a hypothetical self-similar solution to 3D incompressible Euler.
\step
\section{Geometric formulation}\label{sec:geomeqs}

The 3D incompressible Euler equations of fluid dynamics are
\begin{subequations}\label{eq:geomform}
\begin{align}
\mathcal{L}_S t & = 1 \label{transportt}\\
\mathcal{L}_S x^A & = v^A\label{transportx}\\
\mathcal{L}_S \Omega & = 0 \label{transportomega}\\
\dd \Omega & = 0 \label{clom}\\
\delta_{AB}\; \dd v^A \wedge \dd x^B & = \Omega \label{elliptic2} \\
\varepsilon_{ABC}\; \dd v^A \wedge \dd x^B \wedge \dd x^C & = 0 \label{elliptic3}
\end{align}
\end{subequations}
with the physical interpretation\,/\,terminology
\begin{alignat*}{4}
t &\;\;\;\;\;\; && \text{time}\\
x^A &&& \text{position, Cartesian components}\\
v^A &&& \text{velocity, Cartesian components} \displaybreak[0]\\
S &&& \text{material vector field}\\
\text{integral curve of $S$} &&& \text{fluid particle path}\\
\Omega &&& \text{vorticity}
\end{alignat*}

{\bf Preliminaries.}
Equations \eqref{eq:geomform}
are on a 4-dimensional manifold $M^4$, for seven
real valued functions $t$ and $x^A$ and $v^A$,
a vector field $S$, and a 2-form $\Omega$.
For the purpose of this section, they are smooth.
The operator $\mathcal{L}_S$ is the Lie derivative along $S$.
The indices $A,B,C$ each run over $1,2,3$.
By definition $\delta_{AB}$ is equal to $1$
if $A=B$ and $0$ otherwise, and $\varepsilon_{ABC}$ is totally
antisymmetric with $\varepsilon_{123}=1$.
The symmetries of \eqref{eq:geomform} are
discussed in Appendix \ref{app:sym}.

We require that $(t,x^1,x^2,x^3)$ is a diffeomorphism
$M^4 \to \text{(interval)} \times \R^3$,
corresponding to a fluid that fills $\R^3$.
Fluid particles are not allowed to disappear, that is,
each integral curve of $S$ must intersect all level sets of $t$.
Then the flow generated by $S$ maps level sets of $t$ to level sets of $t$ 
by \eqref{transportt}.

The exterior differential calculus in \eqref{eq:geomform}
is not the 4-dimensional one on $M^4$,
but the 3-dimensional one internal to each level set of $t$.
Therefore $\Omega$ is a 2-form on each such level set,
$\dd x^A$ is a 1-form on each such level set, etc.

\mostimportant{Conservation of vorticity \eqref{transportomega}
is here a primitive;
it can be derived
from the more physical velocity-pressure formulation
of incompressible Euler, and is due to Helmholtz and Kelvin.
Conservation of volume
$\mathcal{L}_S (\dd x^1 \wedge \dd x^2 \wedge \dd x^3) = 0$
follows from \eqref{transportx}, \eqref{elliptic3}.
The pressure at unit density is obtained from
$\delta_{AB}\,(\mathcal{L}_S v^A)\, \dd x^B = -\dd \text{(pressure)}$
by integration; the left hand side is closed by
\eqref{transportx}, \eqref{transportomega}, \eqref{elliptic2}.}
\step

{\bf Gauges.} A gauge for \eqref{eq:geomform}
is a choice of independent variables (coordinates).
Different parts of \eqref{eq:geomform} prefer different gauges:
\begin{center}
\begin{tabular}{l|l|l}
\emph{equations} & \emph{kind of problem} & \emph{preferred gauge}
\\
\hline
\eqref{transportt},\eqref{transportx},\eqref{transportomega}
& transport
& Lagrangian
\\
\eqref{clom},\eqref{elliptic2},\eqref{elliptic3}
& elliptic & Eulerian
\end{tabular}
\end{center}

The {\bf Eulerian gauge} promotes $t,x^A$ to independent variables.
In this gauge, the elliptic problem
\eqref{elliptic2}, \eqref{elliptic3} for $v$
is translation invariant.
Equations \eqref{transportt}, \eqref{transportx} imply $S = \p_t + v^A \p_{x^A}$,
which depends on the unknown $v$.

The {\bf Lagrangian gauge}
uses independent variables $t,\lambda^A$, where
the $\lambda^A$ are particle labels, that is,
they are constant along fluid particle paths:
$\mathcal{L}_S \lambda^A = 0$.
Then \eqref{transportt} fixes $S = \p_t$,
which is convenient for the transport problem.
The $\lambda^A$ map each level set of
$t$ diffeomorphically to, by convention, $\R^3$.
One can re-label particles
by acting on $\lambda^A$ with a
time independent diffeomorphism of $\R^3$.

\mostimportant{A common approach uses both gauges:
the Lagrangian gauge for transport,
the Eulerian gauge for the elliptic problem,
with explicit transformations between the two.
Section \ref{sec:locex} pursues another,
consistently Lagrangian approach
to all parts of \eqref{eq:geomform},
including the elliptic problem.}

%%%%%%%%%%%%%%%%%%%%%%%%%%%%%%%%%%%%%%%%%%%%%%%%%%%%%%%%%%%%%%%
\section{Incompressible Euler is a real analytic ODE}\label{sec:locex}
We give a new proof of:
\emph{In the Lagrangian gauge, 3D incompressible Euler is an
autonomous real analytic ODE in Banach space.}

This ODE was already used in \cite{lichtenstein}, \cite{gunther}.
The qualification `real analytic' seems to be due to Serfati \cite{serfati},
with other accounts in \cite{shnirelman}, \cite{inci}.
This section contains a new proof
that is consistently Lagrangian.

\mostimportant{%
In the Lagrangian gauge,
the elliptic problem is not translation invariant, but close to
translation invariant, and the rest
can be controlled by a Neumann series.
Incompressible Euler becomes the ODE \eqref{eq:neumann} below.}

\mostimportant{%
In this paper we always refer to the Lagrangian gauge for local existence.
A solution of \eqref{eq:geomform} is then thought of
as an equivalence class of solutions in the Lagrangian gauge
(solution curves of a real analytic ODE),
two being equivalent iff related by a particle re-labeling.
}

\mostimportant{%
The Lagrangian gauge is naturally interpreted in terms
of particles. In Section \ref{sec:basicargXXXX} there is an equation for $\p_t x$,
but none for $\p_t v$, because we exploit
conservation of vorticity (2 conserved quantities per particle)
and conservation of volume (1 conserved quantity per particle).}

%%%%%%%%%%%%%%%%%%%%%%%%%%%%%%%%%%%%%%%%%%%%%%%%%%%%%%%%%%%%%%%%%%

\subsection{Setup in the Lagrangian gauge}\label{sec:basicargXXXX}
\mostimportant{%
We abbreviate partial derivatives as usual:
if $a,b,c,d$ are used as independent variables,
then $\p_a$ is an abbreviation for
$(\p_a)_{b,c,d}$, the derivative w.r.t.~$a$ at constant $b,c,d$.}

We use Lagrangian coordinates $t,\lambda^A$ as independent variables,
$S = \p_t$.
The vorticity
$\Omega = \Omega_{AB}(\lambda)\;\dd \lambda^A\wedge \dd \lambda^B$
has coefficients independent of $t$ by  \eqref{transportomega},
and is closed by \eqref{clom}. Fix such an $\Omega$.
System \eqref{eq:geomform} reduces to
\begin{subequations}\label{eq:lagformxxx}
\begin{align}
\p_t x^A & = v^A\label{lagode}\\
\delta_{AB}\; \dd v^A \wedge \dd x^B & = \Omega \label{lagell2}\\
\varepsilon_{ABC}\; \dd v^A \wedge \dd x^B \wedge \dd x^C & = 0\label{lagell3}
\end{align}
\end{subequations}

\mostimportant{The unknown is $x = x(t,\lambda)$,
while $v$
is merely a placeholder for the solution to the elliptic problem
\eqref{lagell2}, \eqref{lagell3},
which under suitable assumptions is a map:
given $x$ at time $t$ it determines $v$ at time $t$.
Now \eqref{lagode} is naturally interpreted as an ODE in function space.
System
\eqref{eq:lagformxxx} is invariant under re-labeling,
that is, applying a diffeomorphism of $\R^3$ to $\lambda$.
Below we rewrite the equations, and then this symmetry will
be hidden.
}

Set $x^A = \lambda^A + y^A$, with $y$ the new unknown.
System \eqref{eq:lagformxxx} becomes
\begin{samepage}
\begin{subequations}\label{eq:lagformyyy}
\begin{alignat}{6}
\p_t y^A & = v^A\label{lagodeyyy}\\
\delta_{AB}\; \dd v^A \wedge \dd \lambda^B & =
[\dd v \wedge \dd y]  && +  \Omega \label{lagell2yyy}\\
\varepsilon_{ABC}\; \dd v^A \wedge \dd \lambda^B \wedge \dd \lambda^C & =
[\textstyle\sum_{Q = \lambda,y} \dd v \wedge \dd y \wedge \dd Q] \label{lagell3yyy} && + 0
\end{alignat}
\end{subequations}
where within square brackets $[\;\;]$ we do not keep track of constant coefficients and indices.
\end{samepage}%
Note that \eqref{lagell2yyy}, \eqref{lagell3yyy}
is a linear system for $v$, arranged in the form
$$L_0v = L_y v + (\Omega,0)$$
where $L_0$ and $L_y$ are linear first-order differential operators.
The operator $L_0$ is translation invariant in $\lambda$
and can be inverted
with the Fourier transform.
The
operator $L_y$ has only terms linear and terms quadratic in $y$,
and therefore is small if $y$ is small;
$y$ can be made small initially without
loss of generality by re-labeling.
One can then hope to use a Neumann series,
\begin{equation}\label{eq:neumann}
(\p_ty=)\qquad v = \frac{1}{1-L_0^{-1}L_y} L_0^{-1} (\Omega,0)
\end{equation}
\mostimportant{%
Mere `boundedness' of
$L_0^{-1}L_y$ is due to elliptic regularity:
The loss of a derivative through $L_y$ is regained by $L_0^{-1}$.
It allows one to move some terms that contain $v$ to the right,
without spoiling the analysis.
Elliptic regularity is encoded in assumptions (a.2) and (a.3) below.
}

\mostimportant{%
We are sloppy, especially in (a.3), about the fact that
$L_0$ only returns,
and $L_0^{-1}$ only takes,
pairs $(\omega^2,\omega^3)$ for which $\dd \omega^2 = 0$.
In \eqref{eq:neumann} this is fine,
because $L_0^{-1}$ is only ever applied to $(\Omega,0)$
and to the image of $L_y$.}

%%%%%%%%%%%%%%%%%%%%%%%%%%%%%%%%%%%%%%%%%%%%%%%%%%%%%%
\subsection{Argument}\label{sec:kksiuhikas}
Assumptions:
\begin{itemize}
\itemsep 1pt
\item[(a.1)] Two Banach spaces $\SN$ and $\SDN$
are fixed.
The norms $\NN{\,\cdot\,}$ and $\NDN{\,\cdot\,}$
take functions or tuples of functions
of $\lambda \in \R^3$, depending on context. 
\item[(a.2)]
$\NDN{\p_{\lambda} f} \lesssim \NN{f}$.
\item[(a.3)] 
$\NN{L_0^{-1}f} \lesssim \NDN{f}$.
\item[(a.4)]
$\NDN{fg} \lesssim \NDN{f} \NDN{g}$.
\item[(b)] $\SBN \subset \SN$ is the open subset
of all $f\in\SN$ with $\NDN{\p_{\lambda} f} < \text{(a small constant)}$.
\item[(c)] The fixed vorticity satisfies $\NDN{\Omega} < \infty$.
\end{itemize}
\mostimportant{%
The norm of $\Omega$ is the norm of its coefficients with respect
to the $\dd \lambda^A$. The small constant in (b) is chosen
below, and this choice depends only on 
the structure of the equations and the norms,
but not on $\Omega$.
For each $\Omega$ there is an ODE on $\SBN$.
}

We show that \eqref{eq:neumann}
defines a map
$\SBN \to \SN, y \mapsto v$, with
the uniform bound
\begin{equation} \label{eq:boundvv31}
\NN{v}\,\lesssim\,\NDN{\Omega}
\end{equation}
This map is real analytic, in fact a Taylor series,
and \eqref{lagodeyyy} is an autonomous ODE on $\SBN$
with unique local real analytic solutions $t \mapsto y(t) \in \SBN$.

First note the operator norm
$\NNOP{L_0^{-1}L_y} \lesssim \NDN{\p_{\lambda} y}
+ \NDN{\p_{\lambda} y}^2$
by (a.2), (a.3), (a.4).
Choose the constant in (b) small enough to make
$\NNOP{L_0^{-1}L_y} \leq \tfrac{1}{2}$ for all $y \in \SBN$.
Then the Neumann series converges and we get \eqref{eq:boundvv31}.

To get a Taylor series,
decompose $L_y$ into parts linear
and quadratic in $y$, and then
reorganize the terms produced by the Neumann series
by their homogeneity:
\begin{equation}\label{eq:tyl}
1/(1-L_0^{-1}L_y) \;=\; \textstyle\sum_{n \geq 0} O_n(y,\ldots,y)
\end{equation}
with multilinear operators that for some $R > 0$ satisfy
$$\NNOP{O_n(y_1,\ldots,y_n)} \;\lesssim\; R^{-n}
\NDN{\p_{\lambda}y_1}\cdots \NDN{\p_{\lambda}y_n}$$
with an implicit constant independent of $n$.
Choose the constant in (b) smaller than $\tfrac{1}{2}R$,
to make the Taylor series \eqref{eq:tyl} converge for $y \in \SBN$.

\mostimportant{%
The argument above is abstract.
To ground the discussion,
let $\SN$ and $\SDN$ be subsets of the real-valued
tempered distributions $\mathcal{S}'$ on $\R^3$,
with continuous embedding in
$C^0\subset \mathcal{S}'$, the bounded continuous functions.
Then partial differentiation $\p_{\lambda}: \mathcal{S}' \to \mathcal{S}'$
and pointwise multiplication $C^0 \times C^0 \to \mathcal{S}'$
are defined as usual.
For all
$y \in \SBN$ the map $\lambda \mapsto x = \lambda + y(\lambda)$ is a $C^1$-diffeomorphism of $\R^3$, if the constant in (b) is chosen small enough.

Explicit spaces $\SN,\SDN$ that satisfy
the assumptions are in Appendix \ref{app:norms}.
}

%%%%%%%%%%%%%%%%%%%%%%%%%%%%%%%%%%%%%%%%%%%%%%%%%%%%%%%%%%%%%%%
\subsection{Re-labeling equivalence classes}\label{sec:RELAB}
Re-labeling is a right-action of the group of diffeomorphisms of $\R^3$
on pairs $(x,\Omega)$ in Lagrangian coordinates $t,\lambda$
that maps solutions to solutions, see \eqref{eq:lagformxxx}.
The orbits are the re-labeling equivalence classes.

A diffeomorphism $\lambda = r(\lambda')$
acts by $(x,\Omega) \mapsto (x',\Omega')$
with $x'(t) = r^{\ast}x(t)$ and $\Omega' = r^{\ast}\Omega$.
The first is equivalent to $x'(t,\lambda') = x(t,\lambda)$, or
$$
\lambda' + y'(t,\lambda') = \lambda + y(t,\lambda)$$

\mostimportant{One can narrow the
action to the subgroup of those $r$
that act as continuous and continuously invertible maps
$\SN \to \SN$ and $\SDN \to \SDN$,
on $y$ and $\Omega$ respectively.}

%%%%%%%%%%%%%%%%%%%%%%%%%%%%%%%%%%%%%%%%%%%%%%%%%%%%%%%%%%%%%%%

\subsection{Re-labeling example}\label{relabalgXX}

\mostimportant{%
The construction in Section \ref{sec:kksiuhikas}
stops once $y$ reaches the boundary of $\SBN$,
but before that happens,
one can re-label to push $y$ far back into $\SBN$.
Here we discuss just one such re-labeling, without estimates.}

Fix a re-labeling time $t$.
To make $y' \approx 0$, stipulate
$\lambda' = \lambda + (y \ast m)(\lambda)$
with $m$ a smooth approximation of the identity.
Equivalently $r(\lambda')=\lambda'+c(\lambda')$ with
$c(\lambda') = - (y\ast m)(\lambda' + c(\lambda'))$,
which can be solved by the contraction mapping principle.
 Now
$y' = (y - y\ast m) \circ r$
and
$\Omega' = r^{\ast}\Omega$.

\mostimportant{
We have introduced $m$ to make $r$
smooth, to avoid pull-backs by rough functions.
One can let $m$ approach the identity.
Note that the re-labeling $r$
depends on the unknown $y$, and on $m$.
We will not use such re-labelings.}

%%%%%%%%%%%%%%%%%%%%%%%%%%%%%%%%%%%%%%%%%%%%%%%%%%%%%%%%%%%%%%%
\section{Approximate solution and nearby true solution}\label{sec:appref}

\mostimportant{%
We show how an approximate solution of \eqref{eq:geomform}
can be used as a reference,
to construct a nearby true solution,
with more control (on say the time of existence)
than a direct application of Section \ref{sec:locex}.
There is gauge freedom in how one identifies
the manifold on which the reference lives,
with the manifold on which the new solution
lives. The gauge used in this section
identifies the manifolds by \emph{identifying their fluid paths}.
In Section \ref{secxx4}
we allow for a class of singular references;
Section \ref{sec:hypoloc} will refer back to Section \ref{secxx4}.}

%-------------------------------------------
\subsection{System \eqref{eq:geomform} relative to a reference}\label{secxx1}
Fix a reference
$(M^4,t,X,V,S,\Omega)$
that satisfies
\eqref{eq:geomform}
up to errors $\mathbf{E}_2$, $\mathbf{E}_3$:
\begin{subequations}\label{eqsecxx1}
\begin{align}
\mathcal{L}_S t & = 1 \label{eq:owhkkjdkkkka}\\
\mathcal{L}_S X^A & = V^A \label{eq:XVref}\\
\mathcal{L}_S \Omega & = 0 \label{eq:somxx1}\\
\dd \Omega & = 0\\
\delta_{AB}\; \dd V^A \wedge \dd X^B & = \Omega + \mathbf{E}_2\\
\varepsilon_{ABC}\; \dd V^A \wedge \dd X^B \wedge \dd X^C & = \mathbf{E}_3
\end{align}
\end{subequations}

\mostimportant{By assumption, four equations 
of \eqref{eq:geomform} hold identically.
They comprise
an underdetermined subproblem that can be solved in various ways.
The algebraic manipulations in this section are
intended for errors $\mathbf{E}_2$, $\mathbf{E}_3$ much smaller than $\Omega$.
}

Consider a true solution $(M^4,t,x,v,S,\Omega)$.
As suggested by the notation,
\emph{this solution uses the reference's $M^4$, $t$, $S$
(gauge that identifies fluid paths) as well as $\Omega$.}
Adopt $(t,X)$ as independent variables. Then
$$
S = \p_t + V^A \p_{X^A}
$$
by \eqref{eq:owhkkjdkkkka}, \eqref{eq:XVref}.
The vector field $S$ is fixed by the reference.
Set
\begin{alignat*}{6}
x(t,X) & = & X      &\,+ &\;y(t,X)&\\
v(t,X) & =\;& V(t,X) &\,+ &\;w(t,X)&
\end{alignat*}
with $y$ and $w$ the new unknown.
Equations
\eqref{transportt}, \eqref{transportomega}, \eqref{clom}
hold by construction, and
system \eqref{eq:geomform} is now equivalent to
\begin{subequations}\label{dkhkfhdlffffaa}
\begin{alignat}{6}
\mathcal{L}_S y^A & = w^A \label{ooop1}\\
 \delta_{AB}\; \dd w^A \wedge \dd X^B
& = 
\big[\textstyle\sum_{P=V,w} \dd P \wedge \dd y\big] && -\mathbf{E}_2 \label{ooop2}\\
\varepsilon_{ABC}\; \dd w^A \wedge \dd X^B \wedge \dd X^C 
& =
\big[\textstyle\sum_{P=V,w}\sum_{Q=X,y} \dd P \wedge \dd y \wedge \dd Q\big]
&& -\mathbf{E}_3 \label{ooop3}
\end{alignat}
\end{subequations}
Superficially similar to \eqref{eq:lagformyyy},
this system is not driven by a full-sized $\Omega$,
but by small $\mathbf{E}_2$ and $\mathbf{E}_3$.
Equation \eqref{ooop1}
can be recast as
\begin{equation}\label{eq:odeagain}
\tfrac{\dd}{\dd t}(y(t) \circ \phi(t,s)) = w(t)\circ \phi(t,s)
\end{equation}
for any fixed $s$.
Here $\phi(s,t): \R^3 \to \R^3$ is the flow map of the reference,
$\phi(t,t) = \text{(identity)}$ and
$\tfrac{\dd}{\dd s} \phi(s,t) = V(s) \circ \phi(s,t)$,
see Appendix \ref{app:flowmap} for details.
\step 

\newcommand{\LE}{\ref{sec:locex}}
\mostimportant{%
For local existence we appeal to Section {\LE}, as follows.
Let $x(t),\Omega(t)$ be as in the current section.
For every $s$
let $x_{\LE,s}(t),\Omega_{\LE,s}$
be the solution produced by Section {\LE}
with initial data $x_{\LE,s}(s) = x(s)$
and vorticity $\Omega_{\LE,s} = \Omega(s)$.
Then the correspondence is $x_{\LE,s}(t) = x(t) \circ \phi(t,s)$
and $\Omega_{\LE,s} = \phi(t,s)^{\ast}\Omega(t)$.
To rewrite this in terms of $y$'s,
use $x(t) = \mathbbm{1} + y(t)$
and $x_{\LE,s}(t) = \mathbbm{1} + y_{\LE,s}(t)$ with $\mathbbm{1}: \R^3 \to \R^3$ the identity.
 Then the correspondence is
$y(t) = (\phi(s,t)-\mathbbm{1}) + y_{\LE,s}(t)\circ \phi(s,t)$,
with $\phi$ fixed by the reference.
In particular $y_{\LE,s}(s) = y(s)$.

Basic assumptions about the reference needed for Section {\LE}
are implicit below, e.g.~that $\NDN{\Omega(t)} < \infty$
and $\NN{\phi(s,t)-\mathbbm{1}} < \infty$, and
composition by $\phi(s,t)$ is a bounded linear map on $\SN$.
We assume that the reference is smooth.
To keep Section {\LE} going
we need $y(t) \in \SBN$, 
for which it suffices that $\NN{y(t)}$
is small.
Estimating $\NN{y(t)}$ is the purpose of Section \ref{secxx2}.
}

%%%%%%%%%%%%%%%%%%%%%%%%%%%%%%%

\subsection{Estimate}\label{secxx2}

Equation \eqref{eq:odeagain} implies
\begin{equation}\label{eq:kkkiq}
\big|\tfrac{\dd}{\dd t} \NN{y(t)}\big| \;\leq\; \NN{w(t)} +
\big|\tfrac{\dd}{\dd s}|_{s=t}
\NN{y(t) \circ \phi(t,s)}\big|
\end{equation}
\mostimportant{%
We use the notation
$|\tfrac{\dd}{\dd t} f(t)| = \limsup_{h \to 0} |(f(t+h)-f(t))/h|$.
Inequality \eqref{eq:kkkiq}
is obtained somewhat formally by
setting $a(t,s) = \NN{y(t)\circ \phi(t,s)}$ and using the chain rule
for $\tfrac{\dd}{\dd t} a(t,t)$.
}

In addition to (a) and (b) from Section \ref{sec:basicargXXXX},
 assume:
\newcommand{\SDK}{\text{DK}}
\newcommand{\NDK}[1]{\|#1\|_{\text{DK}}}
\begin{itemize}
\itemsep 0pt
\topsep 0pt
\parsep 0pt
\item[(e.1)] $\SDK$ is a Banach space.
\item[(e.2)] $\NDN{fg} \lesssim \NDK{f}\NDN{g}$.
\item[(e.3)]
$\big|\tfrac{\dd}{\dd s}|_{s=t} \NN{f \circ\textrm{(flow of $U$)}(t,s)}\big|
\;\lesssim\;
\big(|U(0)|
+ \NDK{\p_X U}\big)\, \NN{f}$.
\item[(f)]
 $\NDN{(\mathbf{E}_2,\mathbf{E}_3)(t)}<\infty$ for all $t$.
\end{itemize}
\mostimportant{%
The space $\SDK$ plays a different role than $\SN,\SDN$.
We anticipate that $\SDK$ appears in the
critical assumption (r.1) in Section \ref{secxx4}, and
that one will choose $\SDK$ to suit (r.1).
Estimates of the form (e.3)
for composition by the flow of a vector field $U$
can be obtained from Appendix \ref{app:compest}
if $\SN$ is a weighted Sobolev space;
one does not lose derivatives in $f$.
In (f) the norms of $\mathbf{E}_2$, $\mathbf{E}_3$
are the norms of their coefficients with respect the $\dd X^A$.}

The linear system \eqref{ooop2}, \eqref{ooop3} for $w$ has the form
$$L_0w = L_yw + K_y - (\mathbf{E}_2,\mathbf{E}_3)$$
The operators $L_0$ and $L_y$ are the same as in Section \ref{sec:locex},
except that the independent variables now have a different name,
$X$ instead of $\lambda$. Terms without $w$ go into $K_y$,
they all contain one factor of $\dd y$ at least.
We estimate the map
$\SBN \to \SN, y \mapsto w$ by using a Neumann series
as before,
\begin{equation}\label{eq:uuz}
\NN{w} \;\lesssim\; \NDK{\p_X V} \NN{y} + \NDN{(\mathbf{E}_2,\mathbf{E}_3)}
\end{equation}
where the first term on the right bounds $\NDN{K_y}$.
From \eqref{eq:kkkiq} and \eqref{eq:uuz} we obtain,
uniformly on $\SBN$ and with implicit constant independent of $t$:
\begin{equation}\label{eq:xxxxx}
\big|\tfrac{\dd}{\dd t} \NN{y(t)}\big| \;\lesssim\;
\big(
|V(t,0)| 
+ \NDK{(\p_X V)(t)}\big)\, \NN{y(t)} + \NDN{(\mathbf{E}_2,\mathbf{E}_3)(t)}
\end{equation}

\mostimportant{The map $y \mapsto w$ is time dependent,
through $V$ and $\mathbf{E}_2$ and $\mathbf{E}_3$.
However $L_0$ and $L_y$ do not depend on time,
and neither do the estimates for the Neumann series, they
rely only on $y \in \SBN$.
Therefore the implicit constants in \eqref{eq:uuz}
and \eqref{eq:xxxxx}
are also independent of time, the time dependence is
explicit in the terms on the right hand side.}

%%%%%%%%%%%%%%%%%%%%%%%%%%%%%%%%%%%%%%%%%%%%%%%%%%%%%%%%%%
\subsection{Rescaled estimate}
\newcommand{\resc}{b}
Set $\sigma(\resc): \R^3 \to \R^3, X \mapsto \resc X$.
For every constant $\resc  > 0$,
\begin{multline}\label{eq:dkhkdhdk}
\big|\tfrac{\dd}{\dd t} \NN{y(t) \circ \sigma(\resc )}\big|
\;\lesssim
\big(
\resc^{-1} |V(t,0)|
+ \NDK{(\p_X V)(t) \circ \sigma(\resc )}\big)\, \NN{y(t) \circ \sigma(\resc )}\\ + \resc  \, \NDN{(\mathbf{E}_2,\mathbf{E}_3)(t) \circ \sigma(\resc)}
\end{multline}
provided $y(t) \circ \sigma(\resc ) \in \SBN$.
The implicit constant is independent of $\resc $.

\mostimportant{
Inequality \eqref{eq:dkhkdhdk}
follows from \eqref{eq:xxxxx} by
the scaling symmetry in $X$ of \eqref{dkhkfhdlffffaa}.
The notation
$\mathbf{E}_2(t) \circ \sigma(\resc )$ is short for
$\text{(coefficients of $\mathbf{E}_2$ w.r.t.~the $\dd X$)(t)}
\circ \sigma(\resc )$, and similar for $\mathbf{E}_3$.
We assume that for every $\resc >0$, composition by $\sigma(\resc )$ is
a bounded linear map on $\SN$, $\SDN$, $\SDK$. We only state
this explicitly when we need a quantitative bound.
}

%%%%%%%%%%%%%%%%%%%%%%%%%%%%%%%%%%%%%%%%%%%%%%%%%%%%%%%%%%%%%%
\subsection{Possibly singular reference}\label{secxx4}

\mostimportant{%
Suppose the reference is singular as $t \uparrow 0$.
Under some assumptions we show that
there is a nearby true solution that
is asymptotic to the reference as $t \uparrow 0$.
We do so by solving the equations in the negative time direction,
away from the singularity.
}

Assume that for some constant $A > 0$ the norm satisfies
\begin{itemize}
\item[(q)] $\NN{f} \leq \resc ^{-A} \NN{f \circ \sigma(\resc )}$ for all $\resc  \leq 1$.
\end{itemize}
Assume that the reference is defined for $t$ in some interval
$(-\text{const},0)$,
and that for some constants $B,E,R_1,R_2 > 0$ it satisfies
\begin{itemize}
\item[(r.1)]
 $
|t|^{-B} |V(t,0)|
+ \NDK{(\p_X V)(t) \circ \sigma(|t|^B)} \,\leq\, R_1 |t|^{-1}$.
\item[(r.2)] $\NDN{(\mathbf{E}_2,\mathbf{E}_3)(t) \circ \sigma(|t|^B)} \,\leq\, R_2 |t|^E$.
\end{itemize}

Set
$\mathbf{J}(t) = |t|^{-AB} \NN{y(t) \circ \sigma(|t|^B)}$.
Then
\begin{equation}\label{eq:qoiwele}
\textstyle\limsup_{h \uparrow 0}\;(\mathbf{J}(t+h)-\mathbf{J}(t)) \big/ |h|
\;\lesssim\; R_1 |t|^{-1} \mathbf{J}(t) + R_2 |t|^{B(1-A) + E}
\end{equation}
provided $y(t) \circ \sigma(|t|^B) \in \SBN$.
This follows from
(q), \eqref{eq:dkhkdhdk}, (r.1), (r.2).

\mostimportant{This limits the growth of $\mathbf{J}(t)$
only in the negative time direction:
Given $\mathbf{J}(t)$ with $t < 0$, it yields a bound for
$\mathbf{J}(s)$
with $s < t$.
The implicit constants in \eqref{eq:qoiwele}
in front of the first and second terms,
respectively,
are the same as the implicit constants in \eqref{eq:xxxxx}:
they do not depend on $t$,
and they do not depend on the reference.
In \eqref{eq:qoiwele}
the dependence on the reference
is explicit in $B$, $E$, $R_1$, $R_2$.}

Let $C_1>0$ be the implicit constant in front of the first
term in \eqref{eq:qoiwele}.
Set $\mathbf{K}(t) = |t|^{-C_1R_1} \mathbf{J}(t)$. Then
\begin{equation}\label{eq:oeiuzeizwe}
\textstyle\limsup_{h \uparrow 0}\;(\mathbf{K}(t+h)-\mathbf{K}(t)) \big/ |h|
\;\lesssim\; R_2 |t|^{B(1-A) + E - C_1R_1}
\end{equation}
which is integrable over intervals
$(-\text{const},0)$ iff the exponent is $>-1$, that is
\begin{equation}\label{eq:sjkkkqn}
E \; > \; B(A-1) + C_1R_1 - 1 
\end{equation}
If \eqref{eq:sjkkkqn} holds,
then by a limiting procedure there is a true
solution on some time interval $(-\text{const},0)$
that is asymptotic to the reference as $t \uparrow 0$,
$$\NN{y(t) \circ \sigma(|t|^B)} = \mathcal{O}(|t|^{B+E+1})$$
which includes
$y(t) \circ \sigma(|t|^B) \in \SBN$ for $t$ close to $0$.

\mostimportant{%
The limiting procedure works as follows.
Introduce a family $\{y_s\}$ of solutions parametrized by small $s<0$,
namely set $y_s(s) = 0$ and solve the equations
in the negative time direction. Let $(T(s),s]$
be the domain of $y_s$.
Then there is a common nontrivial time of existence,
$\limsup_{s\uparrow 0} T(s) < 0$, and furthermore $\lim_{s \uparrow 0}y_s$ exists.
To show convergence,
one estimates the distance between pairs
of solutions in the family.
}

%%%%%%%%%%%%%%%%%%%%%%%%%%%%%%%%%%%%%%%%%%%%%%
\section{A hypothetical self-similar flow $\SSS$}\label{sec:sss}

We propose and motivate properties of a hypothetical self-similar
solution $\SSS$ of 3D incompressible Euler.
It would be interesting to know if this solution
exists,
or whether the properties are inconsistent.

\mostimportant{%
If one contemplates consistency,
then this section can be the starting point for 
a numerical search (if a more conceptual approach does not come to mind).
The properties can be used to narrow the numerical search space,
and the section was written expressly with that purpose in mind.
Some properties seem forced,
while others are only motivated as being
the seemingly `simplest' among a set of alternatives,
the idea being that one may want to consider those first.

The perspective taken in this section is that one is looking
for an essentially unique object;
this is spelled out in property \sref{spec:localunique}.
This perspective motivates many of the proposals.

In this section,
we pretend that the hypothetical $\SSS$ exists.
For convenience, we will not endlessly repeat the word hypothetical.
The entire section is understood to be hypothetical.
}

\mostimportant{%
The solution $\SSS$ is self-similar, it develops a singularity,
and it has infinite energy, all of which follows from self-similarity.
Is the singularity peculiar to
this self-similar solution,
or does it encode a mechanism for singularity formation
that extends to finite energy solutions?
A tentative discussion is in Section \ref{sec:hypoloc}.
}

%%%%%%%%%%%%%%%%%%%%%%%%%%%%%%%%%%%%%%%%%%%%%%%%%%%%
\step

This section
has features in common with
an interesting paper by Lathrop \cite{lathrop},
especially with the more mathematical Sections 2, 5, 6, 7 of \cite{lathrop}.

\mostimportant{
Lathrop emphasizes a vector field that he calls $H$,
which is the infinitesimal generator
of the dynamical system \eqref{eq:iiuuw2} below.
Lathrop considers a velocity that goes like
$r^{1-1/\alpha}$ as $r \to \infty$,
in our notation $|x|^{1-1/\KA}$ ($\KA$ is our symbol for $\alpha$)
which is consistent with this section.
The main differences with \cite{lathrop} are as follows.
We allow for
discrete self-similarity, see \sref{spec:selfsim}.
The expansions in and after \sref{spec:iiuu} are new.
The local uniqueness hypothesis \sref{spec:localunique} is new.
It motivates us to narrow the discussion in several places,
say in \sref{spec:conjclass}, and in \sref{spec:XSAY}.
}
\step

An interesting instance of self-similarity
is the general relativity solution found by
Choptuik \cite{choptuik}.

\mostimportant{%
 Choptuik's solution is a pair $(g,\phi)$
with $g$ the metric, $\phi$ a massless scalar field.
The self-similarity is discrete, with a generator that
not only rescales $g$, but also inverts the sign of $\phi$.
That is, the self-similarity is non-trivially coupled to other symmetries.
}

\mostimportant{%
A computer assisted proof of the existence and of the real analyticity
of Choptuik's solution is in \cite{trub}.
This is a case where numerics can be made into a rigorous proof.
It is the solution to an effectively 2-dimensional nonlinear PDE.
A purely conceptual construction is not known at this time.
}

\step

\begin{sproperty}[solution]\label{spec:u1}
The flow $\SSS$ solves incompressible Euler.
Its velocity
$$
(-\infty,0) \times \R^3 \to \R^3,\;\;
(t,x) \mapsto v_{\SSS}(t,x)
$$
with respect to Eulerian coordinates
is smooth, and likely real analytic.
\end{sproperty}

\mostimportant{Eulerian coordinates
are used to simplify the discussion;
the equations for these coordinates are in Appendix \ref{app:euleriangauge}.
It may be that this gauge is not the best
to search for $\SSS$, and that one should for example
use a gauge that `fixes the fluid paths', see below.}

\begin{sproperty}[nontrivial] \label{spec:u2}
The vorticity $\omega_{\SSS} = \nabla \times v_{\SSS}$ does not
vanish identically.
\end{sproperty}

\begin{sproperty}[sublinear growth]\label{sublineargrowth}
$v_{\SSS} = \mathcal{O}(|x|^{1-\delta})$
and $\nabla v_{\SSS} = \mathcal{O}(|x|^{-\delta})$
as $|x|\to \infty$, for fixed $t < 0$.
The constant $\delta > 0$ is independent of $t$.
\end{sproperty}

\mostimportant{
Here and in similar instances below,
`for fixed $t < 0$' is short for
`uniformly on every compact time subinterval of $(-\infty,0)$'.

\sref{sublineargrowth} excludes special
solutions of the form $v = |t|^{-1}
\text{(constant matrix)}\,x$.
They are doubly self-similar,
$v(at,bx) = (b/a) v(t,x)$
for all $a,b \in \R^+$,
and they are the only such solutions.

\sref{sublineargrowth}
guarantees that the flow map (see Appendix \ref{app:flowmap}) exists,
and therefore that Lagrangian coordinates exist.
In fact, \sref{sublineargrowth} implies that the flow map is
asymptotic to the identity:
$\phi^{\SSS}(s,t)x = x + \mathcal{O}(|x|^{1-\delta})$
and its Jacobian is $\mathbbm{1} + \mathcal{O}(|x|^{-\delta})$,
for all fixed $s,t < 0$.
}

\begin{sproperty}[self-similar]\label{spec:selfsim}
There is
a constant $\KA > 0$
and a nontrivial subgroup $Q \subset \R^+$
such that
$v_{\SSS}(qt,\,q^{\KA} x) \;=\; q^{\KA-1}\,v_{\SSS}(t,x)$
for all $q \in Q$.
\end{sproperty}

\mostimportant{
Further symmetries are discussed later.
Most of the discussion uses only \sref{spec:selfsim}.

In a numerical search for $\SSS$, the number $\KA$ is one of the unknowns.

We assume that either $Q = \R^+$, or that $Q$ is isomorphic to $\Z$.
For $Q = \R^+$ the reduced equations are in Appendix \ref{app:css},
but we will not use them.
From now on, all occurrences of $q$ are
understood to be restricted by $q \in Q$.

\sref{spec:selfsim}
implies 
$\phi^{\SSS}(qs,qt) \circ q^{\KA} = q^{\KA} \circ \phi^{\SSS}(s,t)$
where, abusing notation, $q^{\KA}$ denotes the map $x \mapsto q^{\KA}x$.
The vorticity satisfies
$\omega_{\SSS}(qt,q^{\KA} x) = q^{-1} \omega_{\SSS}(t,x)$,
the vorticity 2-form
$(q^{\KA})^{\ast} \Omega^{\SSS}(qt) = q^{2\,\KA-1} \Omega^{\SSS}(t)$.
Also
$(\nabla^k v_{\SSS})(qt,q^{\KA}x) = q^{(1-k)\KA-1} (\nabla^k v_{\SSS})(t,x)$.
}

\begin{sproperty}[flow estimate]\label{hhzzeYY}
If $t < s < 0$ and $|t|^{-\KA}|x| > \text{(big constant)}$ then
\begin{subequations}\label{eq:iiuBOTH}
\begin{align}
\label{eq:iiu11}
|\phi^{\SSS}(s,t)x - x| &
\;\lesssim\; (|t|^{-\KA}|x|)^{-\delta} |x|\\
\label{eq:iiu22}
|\text{(Jacobian of $\phi^{\SSS}(s,t)$ at $x$)} - \mathbbm{1}| &
\;\lesssim\; (|t|^{-\KA}|x|)^{-\delta}
\end{align}
\end{subequations}
\end{sproperty}

\mostimportant{%
Here $\delta$ is the constant in \sref{sublineargrowth}.
Here $A \lesssim B$ is short for $A \leq \text{(big constant)}\,B$.
It is understood that
the big constants do not depend on $s,t$ or $x$.

The assumptions and conclusions in \sref{hhzzeYY}
scale as in \sref{spec:selfsim}, and it therefore suffices to check
\sref{hhzzeYY} for fixed $t < 0$, say $t \approx -1$,
uniformly for all $s \in (t,0)$.
One can use the bound
$|\nabla^k v_{\SSS}|
\lesssim |t|^{\KA\,\delta - 1} \max(|t|^{\KA},|x|)^{1-k-\delta}$
for $k=0$ and $k=1$,
which follows from
\sref{sublineargrowth} and \sref{spec:selfsim}.
}

%--------------------------------------------------------

\begin{sproperty}[almost all particles go to self-similar infinity]\label{uuiew2113}
For almost all fluid paths $t \mapsto x(t)$
one has
$\textstyle\lim_{t \uparrow 0} |t|^{-\KA}|x(t)| = \infty$.
\end{sproperty}

\mostimportant{%
By definition, a fluid path has property $H$ iff
$\textstyle\lim_{t \uparrow 0} |t|^{-\KA}|x(t)| = \infty$.
Property $H$ is equivalent to $\liminf_{t \uparrow 0} |x(t)| > 0$
by \eqref{eq:iiu11}.
Property $H$ is also equivalent to
$\exists t < 0: |t|^{-\KA}|x(t)|
\geq \text{(some threshold)}$, again by \eqref{eq:iiu11}.

Let $a(t)$ be the volume occupied at time $t$ by all particles
that \emph{do not}
have property $H$. Then $a(t)$ is at most the volume of a ball
with radius $\text{(threshold)} |t|^{\KA}$, which goes to zero as
$t \uparrow 0$. Since $a(t)$ is independent of time by incompressibility,
it follows that $a(t)=0$.
}

%--------------------------------------------------------

\newcommand{\AAA}{X}
\begin{sproperty}[limiting flow]\label{h9632Y}
Let $\AAA(r) \subset \R^3$ be the set of points with $|x| > r$.
Set $r(t) = \text{(big constant)} |t|^{\KA}$.
Then
$\lim_{s \uparrow 0} \phi^{\SSS}(s,t)|_{\AAA(r(t))}$
exists,
with $C^1$ convergence on compact subsets,
and is a $C^1$ diffeomorphism $\AAA(r(t)) \to J$ whose image
satisfies $\AAA(2r(t)) \subset J \subset \AAA(r(t)/2)$.
Its Jacobian is uniformly close to $\mathbbm{1}$
and has determinant equal to $1$.
\end{sproperty}

\mostimportant{%
By \sref{spec:selfsim} is suffices to check
this for fixed $t < 0$.
Let $t < s' < s < 0$.
For $C^0$ convergence
use $|\phi^{\SSS}(s,t)x-\phi^{\SSS}(s',t)x|
= |\phi^{\SSS}(s,s')y - y|$ with $y = \phi^{\SSS}(s',t)x$
to apply \eqref{eq:iiu11} and let $s',s \uparrow 0$.
For $C^1$ convergence use \eqref{eq:iiu22} in a similar way,
with the chain rule for $\phi^{\SSS}(s,s')\circ \phi^{\SSS}(s',t)$.

One can derive more properties of
$\lim_{s \uparrow 0}\, \phi^{\SSS}(s,t)$
based on other properties proposed in this section.
Informally speaking, this limit is a degenerate endpoint in
V.I.~Arnold's variational interpretation
of incompressible Euler.
}

%--------------------------------------------------------

\newcommand{\smap}{\Gamma^{\SSS}}
\newcommand{\smapp}{$\Gamma^{\SSS}$ }

\mostimportant{%
To go on, introduce for all $q \in Q$ and $t < 0$ the diffeomorphism
\begin{equation}
\smap_q(t) = q^{-\KA} \circ \phi^{\SSS}(qt,t)\;:\; \R^3 \to \R^3
\end{equation}
By \sref{spec:selfsim} it is a homomorphism in $q$,
and different $t$ yield conjugate maps, explicitly:
\begin{subequations}\label{uausbbaksdq}
\begin{align}
\label{eq:iiuuw2}
\smap_{pq}(t) & = \smap_p(t) \circ \smap_q(t) \qquad p,q \in Q\\
\smap_q(t) & = \phi^{\SSS}(s,t)^{-1} \circ \smap_q(s) \circ \phi^{\SSS}(s,t)
\end{align}
\end{subequations}
Vorticity and volume conservation imply identities internal to level sets of $t$:
\begin{subequations}
\begin{align}
\label{eq:zertzweefe}
\smap_q(t)^{\ast}\Omega^{\SSS}(t) & = q^{1-2\,\KA}\Omega^{\SSS}(t)\\
\label{eq:ooowihebf}
\smap_q(t)^{\ast}(\dd x^1 \wedge \dd x^2 \wedge \dd x^3) & = q^{-3\,\KA}
(\dd x^1 \wedge \dd x^2 \wedge \dd x^3)
\end{align}
\end{subequations}}

\mostimportant{%
To picture $\smap_q$ geometrically,
quotient $(-\infty,0)\times \R^3$
by identifying $(qt,q^{\KA}x) \sim (t,x)$.
The quotient is topologically $S^1 \times \R^3$ if $q \neq 1$,
and it is ruled by fluid paths
that monotonically wind around $S^1$.
Going once around $S^1$ along fluid paths permutes the fluid paths
(pick a transversal hypersurface and look at the return map)
and this permutation is  $\smap_q$, up to conjugacy.
}

\begin{sproperty}[global vorticity bound]\label{SPEC:globalvorticitybound}
$|\omega_{\SSS}| \lesssim \max(|t|^{\KA},|x|)^{-1/\KA}$.
\end{sproperty}

\mostimportant{%
By \sref{spec:selfsim} it suffices to show
$|\omega_{\SSS}| \lesssim \max(1,|x|)^{-1/\KA}$
for fixed $t<0$.
Let
$q \downarrow 0$ in $Q$ in \eqref{eq:zertzweefe}.
For $|x| > \text{(big constant)}$ one has
$|\smap_q(t)x - q^{-\KA}x| \lesssim q^{-\KA} |x|^{1-\delta}$
uniformly for all $q \leq 1$ in $Q$ by \eqref{eq:iiu11},
and a similar bound for the Jacobian using \eqref{eq:iiu22}.
One obtains \sref{SPEC:globalvorticitybound}.
}

\begin{sproperty}[eigenvalues at fixed points]\label{spec:j2382hi}
If $\omega^{\SSS}(t) \neq 0$ at a fixed point of
$\smap_q(t)$,
then the eigenvalues of the linearization satisfy
$\lambda_1 = q^{-1-\KA}$ and $\lambda_2\lambda_3 = q^{1-2\,\KA}$.
\end{sproperty}

\mostimportant{%
For a 3-dim linear map $A$
the identity
$p_{A \wedge A}(\mu) D = -\mu^3 p_A(D/\mu)$
for characteristic polynomials holds, with $D = \det A$.
Here $D = q^{-3\,\KA}$ by \eqref{eq:ooowihebf},
and $p_{A \wedge A}(q^{1-2\,\KA})=0$ by \eqref{eq:zertzweefe}
since $\Omega^{\SSS}(t) \neq 0$. Therefore $p_A(q^{-1-\KA})=0$.
Then use $\lambda_2\lambda_3 = D/\lambda_1$.
}

\begin{sproperty}[classification of \smapp up to conjugacy]\label{spec:conjclass}\rule{0pt}{0pt} These properties hold:
\begin{itemize}
\item[(a)] $\smap_q(t)$ has exactly one fixed point $x_{\ast}(t)$
if $q \neq 1$.
\item[(b)] $\smap_q(t)$
is conjugate in $\text{Diff}^{\infty}(\R^3)$ to a linear map $A_q$.
\item[(c)] $A_q$ is repelling if $q<1$.
\item[(d)] $A_q$ is diagonalizable over $\C$.
\item[(e)] $\omega_{\SSS}(t) \neq 0$ at $x_{\ast}(t)$, for some $t<0$.
\end{itemize}
\end{sproperty}

\mostimportant{%
\sref{spec:conjclass}
collects assumptions about the classification up to conjugacy,
it describes the seemingly simplest scenario.

As suggested by the notation, the fixed point in (a) does not depend on $q$,
and the linear map in (b) does not depend on $t$, by \eqref{uausbbaksdq}.
In (a) there is at least one fixed point since
$\smap_q(t)x \sim q^{-\KA} x$ as $|x|\to \infty$,
and to have exactly one fixed point seems to be  the simplest alternative.
Global linearizability (b) subsumes (a);
note that $q \mapsto A_q$ is a homomorphism.
 Repelling in (c) means
that $|A_q^n\text{(any nonzero vector)}| \to \infty$ as $n \to \infty$.
The non-degeneracy condition (e) seems reasonable,
assuming that there are no symmetries to prevent it.
}

\mostimportant{The classification of $\smap$ up to conjugacy
is a practical matter if
one uses a gauge that `fixes the fluid paths' on the $S^1 \times \R^3$
obtained by identifying $(qt,q^{\KA}x) \sim (t,x)$.
Such a gauge is probably a fine choice,
because then \eqref{eq:zertzweefe}
and vorticity conservation essentially fix $\Omega^{\SSS}$,
but there is an obstruction:
by fixing the fluid paths,
one fixes the conjugacy class for $\smap$.
If one goes along with \sref{spec:conjclass},
then the fluid paths can be fixed up to 
a small number of parameters, namely the eigenvalues in (d),
and there are no other obstructions.
Note that \sref{spec:j2382hi} applies.

With a gauge that does not fix the fluid paths,
such as the Eulerian gauge,
one can in principle ignore the classification up to conjugacy.
}

\begin{sproperty}[all particles but one go to self-similar infinity]\label{spec:9u39hh93}
The map $t \mapsto x_{\ast}(t)$ from
 \sref{spec:conjclass}
is a fluid path and satisfies $|t|^{-\KA}|x_{\ast}(t)| \lesssim 1$ uniformly
in time.
All other fluid paths $t \mapsto x(t)$ satisfy
$\lim_{t \uparrow 0} |t|^{-\KA}|x(t)| = \infty$.
\end{sproperty}

\mostimportant{%
This is a consequence of \sref{spec:conjclass}.
The properties of $t \mapsto x_{\ast}(t)$ follow
from \eqref{uausbbaksdq}
and from the fixed point condition
$x_{\ast}(qt) = q^{\KA} x_{\ast}(t)$.
For all other fluid particle paths,
(c)
implies that
$\smap_q(t)x(t) = q^{-\KA} x(qt)$
eventually leaves, as $q \downarrow 0$ in $Q$ at fixed $t$,
every compact subset of $\R^3$, which gives
$\lim_{t \uparrow 0} |t|^{-\KA}|x(t)| = \infty$.
This strengthens 
\sref{uuiew2113}.}

\begin{sproperty}[lower bound for $\KA$]\label{xxxuu12}
 $\KA \,>\, \tfrac{1}{2}$.
\end{sproperty}

\mostimportant{%
By \sref{spec:j2382hi} and \sref{spec:conjclass}.
If $q<1$ then the eigenvalues of $A_q$ must satisfy
$|\lambda_1| > 1$, $|\lambda_2| > 1$, $|\lambda_3| > 1$
by (c), (d). Note that (e) holds for all $t$ if it holds for one $t$ by
vorticity conservation, since $t \mapsto x_{\ast}(t)$ is a fluid path.
Then $|\lambda_2\lambda_3| > 1$ and \sref{spec:j2382hi}
imply $1-2\,\KA < 0$.

Incidentally, $\tfrac{1}{2}$ is the only
value compatible with Navier-Stokes diffusion.
A self-similar
solution of incompressible Euler with $\KA > \tfrac{1}{2}$
 is in some sense
an approximate solution of Navier-Stokes at early times
$t \to -\infty$, but the relative strength of diffusion
grows as $t \uparrow 0$.
}

\mostimportant{%
To go on, we assume that $\SSS$ is in some sense
`regular at self-similar infinity'.
The ansatz \sref{spec:iiuu} accommodates
elliptic modes \sref{spec:elliptic}
and transport modes \sref{spec:transport},
and is `closed under multiplication'
in view of the nonlinear nature of incompressible Euler.
Other than that, the ansatz is intended to be minimal.
}

\begin{sproperty}[asymptotic expansion]\label{spec:iiuu}
Set $a_{mn}^{\alpha}(z) = z^{\alpha m+n}$.
Then asymptotically as $|x|^{-1/\KA} |t| \to 0$, with a disclaimer below,
\begin{subequations}\label{eq:llzzuuBOTH}
\begin{alignat}{5}
\label{eq:llzzuu}
v_{\SSS}(t,x) \; & \sim\;\;& |x|/|t|
&\textstyle\sum_{m,n \geq 0}
c[v_{\SSS}]_{mn} (t,x)\, a_{mn}^{\KA}(|x|^{-1/\KA}|t|)\\
\label{eq:llzzuuvo}
\omega_{\SSS}(t,x) \; & \sim\;\;& 1/|t|
&\textstyle\sum_{m,n \geq 0}
c[\omega_{\SSS}]_{mn} (t,x)\, a_{mn}^{\KA}(|x|^{-1/\KA}|t|)
\end{alignat}
\end{subequations}
The coefficient functions
$c[\,\cdot\,]_{mn}: (-\infty,0) \times
(\R^3 \setminus \{0\}) \to \R^3$
are doubly periodic:
$c(qt,x) = c(t,q^{\KA} x) = c(t,x)$ for all $q \in Q$.
\end{sproperty}

\mostimportant{%
Doubly periodic functions
live on a compact space, $S^2$ or $S^1\times S^1\times S^2$ 
depending on $Q$. In particular, doubly periodic implies bounded.
The dimensionless
$|t|\p_t$ and $|x| \nabla$ map doubly periodic to doubly periodic.
%----
Note that \eqref{eq:llzzuuBOTH} is
consistent with \sref{spec:selfsim}, and is
 both
an asymptotic expansion for fixed $t<0$ as $|x| \to \infty$,
and 
for fixed $x \neq 0$ as $t \uparrow 0$.
}

\mostimportant{%
Disclaimer: The $a_{mn}^{\KA}$ are linearly dependent
at rational $\KA$,
which is unacceptable in \eqref{eq:llzzuuBOTH}.
One can replace the $a_{mn}^{\alpha}$ by a linearly independent
but otherwise equivalent family, 
see for example Appendix \ref{app:specialfunc}.
To simplify the discussion, we keep using the $a_{mn}^{\alpha}$.
At rational $\KA$, one can make some sense of this
by letting the $c[\,\cdot\,]_{mn}$
be meromorphic functions of $\alpha$ with rational poles
that cancel upon summation,
as in
$((\alpha-1)^{-1} a^{\alpha}_{10}(z)
  - (\alpha-1)^{-1} a^{\alpha}_{01}(z))|_{\alpha = \KA}$.
Similarly for other coefficients and expansions below.

\sref{spec:iiuu} yields formally
a compactification of the problem:
there is only a discrete set of coefficients $c[\,\cdot\,]_{mn}$
even though the domain of $\SSS$ is unbounded.
It is conceivable that
the asymptotic expansions converge and are equal to $v_{\SSS}$, $\omega_{\SSS}$
 for small 
$|x|^{-1/\KA}|t|$.

In \sref{spec:iiuu} there are two indices $m,n$ for just
one geometric dimension, namely $|x|^{-1/\KA}|t|$.
If correct, then this may rule out
certain numerical approaches to finding $\SSS$.
A numerical search that uses \eqref{eq:llzzuuBOTH},
or is at least aware of it, may have an advantage.

Complications for a numerical search that uses \eqref{eq:llzzuuBOTH}:
the $a_{mn}^{\KA}$ depend on the unknown $\KA$;
one may have to use something like
Appendix \ref{app:specialfunc};
two indices $m,n$ for just
one geometric dimension.
If $Q$ is discrete, then there are two more Fourier indices.
Separate expansions are needed in the bulk, away from self-similar infinity.
}

\begin{sproperty}[sublinear growth for coefficients]\label{spec:triang}
$c[v_{\SSS}]_{00}=0$ and $c[\omega_{\SSS}]_{00}=0$.
\end{sproperty}

\mostimportant{%
By \sref{sublineargrowth}.

Plugging \eqref{eq:llzzuuBOTH}
into incompressible Euler \eqref{rtriziugdu}
yields a problem for the coefficients $c_{mn}$.
We now discuss just this problem.
It is not by itself enough to determine all the $c_{mn}$:
there are `constants of integration'
that can presumably only be gotten by calculating $\SSS$ in the bulk.

This problem has a useful symmetry:
for all $\mu,\nu > 0$,
replacing all $c_{mn}$ by $\mu^m \nu^n c_{mn}$
both for $v_{\SSS}$ and $\omega_{\SSS}$ maps solutions to solutions.
That is, linear terms $c_{mn}$ and quadratic terms
$c_{m'n'}c_{m''n''}$ only appear together in an equation
when $m = m'+m''$ and $n=n'+n''$;
we have suppressed derivatives etc,
and each coefficient belongs to either $v_{\SSS}$ or $\omega_{\SSS}$.

Clearly $m',m'' \leq m$ and $n',n'' \leq n$,
however after using \sref{spec:triang} to drop terms
we also get $(m',n') \neq (m,n)$ and $(m'',n'') \neq (m,n)$.
This yields a recursion
where step $(m,n)$
has the form $A c_{mn} = B$,
with an expression
$B$ that is quadratic in coefficients constructed `earlier',
and with a linear operator $A$ that depends
on $m, n, \KA$, but not on say any of the coefficients.

The linear operator $A$ has nontrivial kernel when
 $(m,n)$ has the form $(m,0)$ or $(0,1)$.
In these cases $B = 0$,
that is, these steps of the recursion are governed by the linear equations
$\p_t \omega = 0$, $\nabla \cdot\omega = 0$, $\nabla \times v = \omega$,
$\nabla \cdot v = 0$ with respectively
\begin{alignat*}{9}
v & =\, &c[v_{\SSS}]_{m0}(t,x)&|x|^{1-m}&& |t|^{\KA\,m - 1}
&\qquad\qquad
v & =\, &c[v_{\SSS}]_{01}(t,x)&|x|^{1-1/\KA}\\
\omega & = &c[\omega_{\SSS}]_{m0}(t,x) &|x|^{-m} &&|t|^{\KA\,m-1}
&
\omega & = &c[\omega_{\SSS}]_{01}(t,x)&|x|^{-1/\KA}
\end{alignat*}
and this notation is not used elsewhere.
The kernels are in \sref{spec:elliptic} and \sref{spec:transport},
they are the `constants of integration' alluded to above.
The operator $A$ is invertible for all other $(m,n)$
if $\KA$ is irrational, for rational $\KA$
see the disclaimer after \sref{spec:iiuu}.}

\FloatBarrier
\begin{sproperty}[decaying elliptic modes]\label{spec:elliptic}
$c[\omega_{\SSS}]_{m0} = 0$,
and the $n=0$ subexpansion of \eqref{eq:llzzuu} is equal to
\begin{equation}\label{eq:kksw12}
(1/|t|)\,\textstyle\sum_{\ell \geq 1}\sum_{-\ell \leq k \leq \ell}
E_{\ell k}^{\SSS}(t)\,|t|^{\KA(\ell+3)}\,
\nabla(|x|^{-\ell-1}Y_{\ell k})
\end{equation}
for some scalar-valued coefficients with $E(qt)=E(t)$.
\end{sproperty}

\mostimportant{%
The spherical harmonics $Y_{\ell k}$ are functions of $x/|x|$
with $\Delta_{S^2} Y_{\ell k} = -\ell(\ell+1)Y_{\ell k}$.
Expression \eqref{eq:kksw12} is for every $t$
the gradient of the multipole expansion of
a general decaying harmonic function in $x$.
There is no $\ell = 0$ term by incompressibility:
the flux of $v_{\SSS}$ through spheres of constant $|x|$
must vanish identically.
For each $\ell \geq 1$,
the $2\ell + 1$ functions $E_{\ell k}^{\SSS}$
determine $c[v_{\SSS}]_{\ell+3,0}$,
while $c[v_{\SSS}]_{m0} = 0$ for $m \leq 3$.
}

\mostimportant{%
Growing modes $\nabla(|x|^{\ell} Y_{\ell k})$
with $\ell \geq 2$
are suppressed by \sref{sublineargrowth},
but not the more subtle $\ell = 1$
which goes into $c[v_{\SSS}]_{10}$.
However, the symmetry `shake' in Appendix \ref{app:sym}, that takes
$v(t,x)$ to $v(t,x-\beta(t)) + \frac{\dd}{\dd t}\beta(t)$,
maps self-similar solutions to self-similar solutions if
$\beta(qt) = q^{\KA} \beta(t)$.
Among the solutions related by this symmetry,
only one has no $\ell = 1$ growing mode,
and \sref{spec:elliptic} selects that one.
Thus `shake' is fixed, and the absence of growing elliptic modes
is a full set of `elliptic boundary conditions at infinity'
for $v_{\SSS}$.

These vanish: $c_{00}[v_{\SSS}]$ by \sref{sublineargrowth},
$c_{10}[v_{\SSS}]$ by `shake',
$c_{20}[v_{\SSS}]$ automatically,
$c_{30}[v_{\SSS}]$ by incompressibility.
}

\begin{sproperty}[transport modes]\label{spec:transport}
$c[v_{\SSS}]_{01}(t,x) = T^{\SSS}(x)$
for a vector-valued $T^{\SSS}$ that does not vanish identically,
$T^{\SSS}(q^{\KA} x) = T^{\SSS}(x)$
and $\nabla \cdot(|x|^{1-1/\KA}T^{\SSS})=0$.
\end{sproperty}

\mostimportant{
Then $c[v_{\SSS}]_{01}$ is
 the leading term of the expansion
assuming only $\KA > \tfrac{1}{4}$, see \sref{xxxuu12}.
}

\begin{figure}[h]
\centering
\captionsetup{justification=centering}
\input{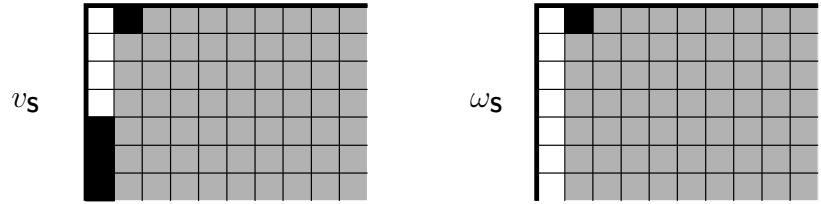}
\caption*{%
\footnotesize
Coefficients $c[v_{\SSS}]_{mn}$
and $c[\omega_{\SSS}]_{mn}$
with $m\geq 0$ vertical, $n\geq 0$ horizontal.\\
White boxes vanish.
Black boxes are elliptic and transport modes.
}
\end{figure}

\FloatBarrier

\begin{sproperty}[global velocity bound]\label{spec:globalvelocitybound}
$|\nabla^k v_{\SSS}| \lesssim \max(|t|^{\KA},|x|)^{1-1/\KA-k}$
where the implicit constant depends on $k=0,1,2,\ldots$
\end{sproperty}

\mostimportant{
By \sref{spec:selfsim} it suffices to
show \sref{spec:globalvelocitybound} for
fixed $t < 0$, which holds by \sref{spec:iiuu}-\sref{spec:transport}
and \sref{xxxuu12}.

Note: If $\KA > 1$ then $v_{\SSS}$ is bounded near
$(t,x)=(0,0)$, by \sref{spec:globalvelocitybound} with $k=0$.
If $\KA < 1$ then $v_{\SSS}$ is unbounded near $(0,0)$,
 by  \sref{spec:selfsim}.
We take no stance on whether $\KA>1$ or not.
}

\begin{sproperty}[final velocity]\label{spec:finalvelocity}
$\lim_{t \uparrow 0} v_{\SSS}(t,x) = |x|^{1-1/\KA} T^{\SSS}(x)$ for
all $x \neq 0$.
\end{sproperty}

\begin{sproperty}[finite energy in bounded domains]
The energy in the ball with radius $r > 0$
stays bounded,
$\limsup_{t \uparrow 0} \int_{|x| < r} \dd^3 x\,  |v_{\SSS}(t,x)|^2 < \infty$.
\end{sproperty}

\mostimportant{%
The integral is $\lesssim \int_{|x|<r}\dd^3x\, \max(|t|^{\KA},|x|)^{2-2/\KA}$.
Split it into respectively $|x| < |t|^{\KA}$ and $|t|^{\KA} < |x| < r$.
One obtains $\lesssim |t|^{B\,\KA} + r^B/B$
provided $B = 5-2/\KA$ is positive.
}

\begin{sproperty}[local uniqueness]\label{spec:localunique}
Informally,
$\SSS$ is uniquely characterized,
locally in function space, by the properties
\sref{spec:u1}, \sref{spec:u2},
\sref{sublineargrowth}, \sref{spec:selfsim}.
Local uniqueness is understood modulo
the symmetries in Appendix \ref{app:sym}.
\end{sproperty}

\mostimportant{%
This makes searching for $\SSS$ a more attractive proposition.
\sref{spec:localunique} proposes that $\SSS$
is not part of a continuous family of solutions.
If there is a discrete family,
then there may still be one member that stands out as the most basic,
say the `most stable'.

\sref{spec:localunique} can be informally
motivated by counting degrees of freedom.
We discuss the case when $Q\subset \R^+$ is discrete,
the case $Q=\R^+$ is similar.
Let $q<1$ be the generator of $Q$.
Consider the fundamental region $[-1,-q] \times \R^3$.
Sublinear growth \sref{sublineargrowth} only suppresses the
growing elliptic modes,
and yields a full set of `elliptic boundary conditions at infinity'
for $v_{\SSS}$;
see also remarks after \sref{spec:elliptic}.
Initial data at $t = -1$ maps,
assuming suitable existence and uniqueness,
to data at $t = -q$, but
only a discrete family of initial data
satisfies the self-similarity condition \sref{spec:selfsim}.
This counting has a flaw: \sref{spec:selfsim}
refers to $\KA$ and the generator $q$,
two more unknowns that seem to spoil the counting.
The remedy are two residual symmetries,
scaling $x$
and scaling $t$,
which is an $\R^+ \times \R^+$ symmetry modulo the discrete $Q$.
Gauge fixing these two symmetries restores the counting.

With \sref{spec:localunique} it would seem that
$\SSS$ should have as much symmetry as possible.
However symmetry is not linearly ordered,
and it is not clear what `possible' means.
One should study symmetry systematically.
A very limited discussion is in
Appendix \ref{app:somestabcand}.
The symmetries proposed in \sref{spec:selfsim}
and \sref{spec:XSAY} are
among those in Appendix \ref{app:somestabcand}.
Prejudicing a numerical search to a more symmetric class makes
things easier resource-wise, but is more dangerous.
}

\begin{sproperty}[symmetry]\label{spec:XSAY}
Let $R_{\varphi} \in \SO(3)$ be rotation around
a fixed axis by angle $\varphi$.
Then the symmetries of $\SSS$ are a subgroup
of $\R^+ \times \{\pm 1\} \times (\R/2\pi \Z)$,
where each element $(p,\sigma,\varphi)$ is the symmetry\; $v_{\SSS}(pt,p^{\KA} \sigma R_{\varphi}x) = p^{\KA-1} \sigma R_{\varphi} v_{\SSS}(t,x)$.
\end{sproperty}

\mostimportant{%
\sref{spec:XSAY} proposes that
a numerical search should not assume symmetries beyond those given.
Such a search does not rule out other symmetries,
but does not use them either.

\sref{spec:XSAY} takes no stance on which subgroup the symmetries of $\SSS$
are. Some possibilities:
\begin{itemize}
\item[(i0)] all elements
\item[(i1)] all elements with $\mu \log p = \varphi$, for a constant $\mu \neq 0$
\item[(i2)] all elements with $p = (p_1)^n$ and $n \in \Z$, for a constant $p_1 > 1$
\item[(i3)] like (i2), but with the additional condition $(-1)^n = \sigma$
\end{itemize}

They
 are all consistent with \sref{spec:selfsim}, explicitly $Q$ is
given by: (i0) $Q = \R^+$, (i1) $ Q= (e^{2\pi k/\mu})_{k \in \Z}$,
(i2) $Q=((p_1)^k)_{k \in \Z}$, (i3) $Q=((p_1)^{2k})_{k \in \Z}$.

The subgroup (i0) has continuous dimension 2,
while (i1), (i2), (i3) have dimension 1.

In (i0), (i1), (i2) there is no
restriction on $\sigma$. One obtains other subgroups by putting $\sigma = 1$.
Assuming the $\{\pm 1\}$ symmetry in a numerical search
yields a small gain at best.}

%%%%%%%%%%%%%%%%%%%%%%%%%%%%%%%%%%%%%%%%%%%%%%%%%%%%%%%%%%%%%%%
\section{A tentative
strategy to localize $\SSS$}\label{sec:hypoloc}

This entire section is under the
assumption that $\SSS$ in Section \ref{sec:sss} exists,
it is self-similar, it develops a singularity, and it has infinite energy.
Does it encode a
mechanism for singularity
formation that extends to finite energy solutions?
This section has no answer, but only a few remarks.

We sketch a strategy that modifies $\SSS$ to a
finite energy flow $\SSS'$
that does not itself solve incompressible Euler,
but is a good enough approximate solution to
invoke Section \ref{secxx4},
which (by solving incompressible Euler backwards)
produces a solution $\SSS''$ that has finite energy and
develops a singularity.

We do not know if this strategy can be implemented.
In step \stepref{underdetproblem} we 
encounter a nonlinear system of equations that looks underdetermined, but
we do not know if there are obstructions to solving it.
The discussion is very informal.

\newcommand{\restrS}{\SSS|_{\text{ball}}}

\newcommand{\EXT}{\text{\bf EXT}}
\newcommand{\LLL}{\text{\bf L}}

\begin{sstep}\label{xosnfjlf} Let $\restrS$
be the restriction of $\SSS$ to a ball in $\R^3$.
\end{sstep}

\mostimportant{%
Along the spherical boundary
and as $t \uparrow 0$,
the asymptotic expansion \eqref{eq:llzzuu} for $\SSS$ applies.}

\begin{sstep}
In the complement of the ball, introduce a formal series flow
\begin{equation}\label{eq:llzzuu333}
v_{\LLL}(t,x) \; = \; (|x|/|t|)
\textstyle\sum_{m,n \geq 0}
c[v_{\LLL}]_{mn} (t,x)\, a_{mn}^{\KA}(|x|^{-1/\KA}|t|)
\end{equation}
Require $c(qt,x) = c(t,x)$ for all $q \in Q$,
\emph{but do not require periodicity in $x$.}
\end{sstep}

\mostimportant{%
Without periodicity in $x$, one could absorb
all powers of $|x|$ in the coefficients, but keeping them
makes it easier to compare with
\eqref{eq:llzzuu}. The formal series \eqref{eq:llzzuu333} is an expansion at $t = 0$.
}

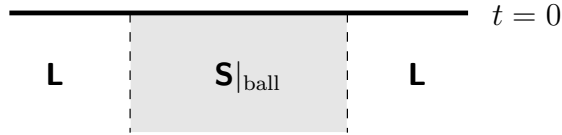
\begin{figure}[h]
\centering
\input{lok.pstex_t}
\end{figure}

\begin{sstep}\label{sstep:elliptic}
Set $c[v_{\LLL}]_{m0} = c[v_{\SSS}]_{m0}$.
\end{sstep}

\mostimportant{%
The corresponding terms in \eqref{eq:llzzuu333} decay like
$\mathcal{O}(|x|^{-3})$,
consistent with finite energy,
because $c[v_{\LLL}]_{m0}=0$ when $m \leq 3$,
see the discussion after \sref{spec:elliptic}.
In particular, $c[v_{\LLL}]_{00}=0$.}

\begin{sstep}
Data for $\LLL$ along $t = 0$ is a function
$c[v_{\LLL}]_{01}(t,x) = T^{\LLL}(x)$. It must
\begin{itemize}
\itemsep 0pt
\item be compatible with incompressibility $\nabla \cdot v_{\LLL} = 0$,
\item match $c[v_{\SSS}]_{01}(t,x) = T^{\SSS}(x)$ smoothly
at the spherical boundary,
\item decay, with its derivatives, sufficiently fast as $|x| \to \infty$.
\end{itemize}
Other than that, it can be chosen freely.
\end{sstep}

\begin{sstep}
Require that the formal series
$\LLL$ solves incompressible Euler,
and that the orthogonal components of
$v_{\LLL}$ and of the asymptotic series of $v_{\SSS}$
match continuously at the spherical boundary.
Now $T^{\LLL}$ determines $v_{\LLL}$
by a recursion.
\end{sstep}

\mostimportant{%
One can only prescribe the orthogonal velocity component
at the boundary.

The choice in \stepref{sstep:elliptic}
is compatible with both incompressible Euler and the boundary condition.
It seems that \stepref{sstep:elliptic} is forced,
because $T^{\LLL}$ does not affect the $c[v_{\LLL}]_{m0}$ components.
}

\begin{sstep}\label{underdetproblem}
Let $K > 0$ be a big integer.
\emph{Can one tweak $T^{\LLL}$ so that for all $|m| \leq K$ and $|n| \leq K$,
the coefficients $c[v_{\LLL}]_{mn}$ and
$c[v_{\SSS}]_{mn}$ match at the spherical boundary,
with $K$ continuous derivatives?}
\end{sstep}

\mostimportant{
This problem seems to be underdetermined.
For example if $Q = \R^+$,
then $T^{\LLL}$ is a function of three independent variables
(exterior of a ball in $\R^3$),
but the matching conditions concern only a finite number
of functions of two independent variables (sphere in $\R^3$).}

\begin{sstep}
Assuming \stepref{underdetproblem} has a solution,
set
$\SSS' = \text{smoothing}(\restrS,\text{truncation}(\LLL))$.
The truncation of $\LLL$ is the finite sum
over $|m|,|n| \leq K$.
It does not match $\restrS$ at the spherical boundary,
not even continuously,
but the mismatch decays quickly as $t \uparrow 0$
given \stepref{underdetproblem}.
To correct the mismatch, apply a smoothing near the boundary,
chosen so as to not destroy the $t \uparrow 0$ asymptotics.
\end{sstep}

\begin{sstep}
The flow $\SSS'$ is only an approximate solution of incompressible Euler,
with an error that decays quickly as $t \uparrow 0$.
One can now invoke Section \ref{secxx4}
(which solves the equations in the negative time direction)
to get a true solution $\SSS''$ on a small
time interval, with a singularity as $t \uparrow 0$.
\end{sstep}

\mostimportant{%
Disclaimer:
It is clear how
Section \ref{secxx4},
is intended to be applied,
but we have not checked all the assumptions there.
Certainly, $K$ must be big enough
in view of condition  \eqref{eq:sjkkkqn}.

Disclaimer: Section \ref{sec:appref} requires that references satisfy four
equations of \eqref{eq:geomform} identically,
see \eqref{eqsecxx1}, and one must adjust
$\SSS'$ accordingly. This should only be a technical matter.

Note that $\SSS'$ coincides with $\SSS$ for small $|x|$,
but $\SSS''$ does not.

We have used the Eulerian gauge to keep the discussion short.
For actual calculations
it may be easier to first put Lagrangian coordinates on
$\SSS$.
Instead of restricting $\SSS$ to a ball in Eulerian coordinates,
one uses a ball in Lagrangian coordinates.
The new boundary
is `characteristic' (ruled by fluid paths)
which may be helpful.

The strategy above, if it works at all,
only produces rather unnatural examples.
It would be interesting to study
more generic classes of finite energy solutions that are
singular like $\SSS$.
}

\begin{comment}

\begin{align*}
\NN{f} & =
  \|\langle x \rangle \p f\|_{L^2}
+ \|\langle x \rangle \p^2 f\|_{L^2}
+ \|\langle x \rangle \p^3 f\|_{L^2}\\
\NDN{f} & =
  \|\langle x \rangle f\|_{L^2}
+ \|\langle x \rangle \p f\|_{L^2}
+ \|\langle x \rangle \p^2 f\|_{L^2}\\
\NDK{f} & = 
  \| f\|_{\infty} + \|\p f \|_{\infty} + \|\p^2 f\|_{\infty}
\end{align*}
\end{comment}
%%%%%%%%%%%%%%%%%%%%%%%%%%%%%%%%%%%%%%%%%%%%%%%%%%%%%%%%%%%%%%%

\section*{Acknowledgments}
This paper is based on work at ETH Zurich and IAS Princeton.
I thank Eugene Trubowitz for many
useful discussions concerning this paper.
Our joint \cite{trub} was one of the motivations for this paper.
I enjoyed support from
The Giorgio and Elena Petronio Fellowship Fund,
and the (US) National Science Foundation.

\mostimportant{This material is based upon work supported by the National Science Foundation under agreement No.~DMS-1128155. Any opinions, findings and conclusions or recommendations expressed in this material are those of the author(s) and do not necessarily reflect the views of the National Science Foundation.}

%%%%%%%%%%%%%%%%%%%%%%%%%%%%%%%%%%%%%%%%%%%%%%%%%%%%%%%%%%%%%%

\mostimportant{%

}

%%%%%%%%%%%%%%%%%%%%%%%%%%%%%%%%%%%%%%%%%%%%%%%%%%%%%%%%%%%%%%%
\appendix

\section{Symmetries of the equations \eqref{eq:geomform}}\label{app:sym}
\FloatBarrier
\begin{figure}
\centering
\begin{tabular}{r|c|c|c|c|c|l}
  & new $t$ & new $x^A$ & new $v^A$ & new $S$ & new $\Omega$ & \\
\hline
$\text{ttrans}_\alpha$ & $t + \alpha$ & $x^A$ & $v^A$ & $S$ & $\Omega$ &\\
$\text{xtrans}_\alpha$ & $t$ & $x^A + \alpha^A$ & $v^A$ & $S$ & $\Omega$ & \\
$\text{rot}_\alpha$ & $t$ & ${\alpha^A}_B x^B$ & ${\alpha^A}_B v^B$ & $S$ & $\Omega$ & 
$\alpha \in \text{O(3)}$\\
$\text{boost}_\alpha$ & $t$ & $x^A + \alpha^A t$ & $v^A + \alpha^A$ & $S$ & $\Omega$ \\
$\text{tscale}_\alpha$ & $\alpha t$ & $x^A$ & $v^A/\alpha$ & $S/\alpha$ & $\Omega/\alpha$ & $\alpha > 0$\\
$\text{xscale}_\alpha$ & $t$ & $\alpha x^A$ & $\alpha v^A$ & $S$ & $\alpha^2 \Omega$ & $\alpha > 0$\\
timeinv & $-t$ & $x^A$ & $-v^A$ & $-S$ & $-\Omega$\\
$\text{shake}_\beta$ & $t$ & $x^A + \beta^A$ & $v^A + \mathcal{L}_S \beta^A$ & $S$ & $\Omega$ & $\dd \beta^A = 0$\\
$\text{diff}_\varphi$ & $\varphi^{\ast}t$ & $\varphi^{\ast}x^A$ &
 $\varphi^{\ast} v^A$ & $\varphi^{\ast}S$ & $\varphi^{\ast}\Omega$ &
$\varphi \in \text{Diff}(M)$
\end{tabular}
\caption*{All $\alpha$ are constants, $\beta$ is a function of time,
$\varphi$ is a diffeomorphism.}
\end{figure}

The symmetries of \eqref{eq:geomform} in the table
act on the set of tuples $(t,x^A,v^B,S,\Omega)$
on a fixed manifold $M$,
and map solution tuples to solution tuples.
The Galilean group is given by
{\bf ttrans}, {\bf xtrans}, {\bf rot}, {\bf boost}.
There are two scaling symmetries
{\bf tscale} and {\bf sscale}
because time and space units
are independent in incompressible Euler,
unlike 
say Navier-Stokes which couples the two units via the
diffusion constant.
By {\bf timeinv} the equations can be solved in either time direction.
The symmetry $\text{\bf shake}$ generalizes xtrans and boost,
and it invalidates a naive uniqueness statement,
since the $\beta^A$ are arbitrary functions of time $t$.
Finally, {\bf diff}
makes explicit the geometric nature of the formulation.
%\FloatBarrier

%%%%%%%%%%%%%%%%%%%%%%%%%%%%%%%%%%%%%%%%%%%%%%%%%%%%%%%%%%%%%%%%%%%%%%%%%%%%%
\section{System \eqref{eq:geomform}
in the Eulerian gauge}\label{app:euleriangauge}
With independent variables $t,x^A$ system \eqref{eq:geomform} reduces to
\begin{subequations}\label{rtriziugdu}
\begin{align}
\label{eq:vcsv}
(\p_t + v \cdot \nabla)\,\omega^A & =
(\omega \cdot \nabla)\,v^A - \omega^A\,(\nabla \cdot v)\\
\nabla \cdot \omega & = 0\\
\nabla \times v & = \omega\\
\nabla \cdot v & = 0
\end{align}
\end{subequations}

\mostimportant{%
They correspond to
\eqref{transportomega}, \eqref{clom}, \eqref{elliptic2}, \eqref{elliptic3}
after substituting
$S = \p_t + v^A \p_A$,
which is equivalent to \eqref{transportt}, \eqref{transportx}.
We have defined $\nabla = (\p_1,\p_2,\p_3)$
and
$\Omega = \tfrac{1}{2}\varepsilon_{ABC}\,\omega^A\, \dd x^B \wedge \dd x^C$.

Equations \eqref{transportomega} and \eqref{eq:vcsv}
are explicitly related by (using $\mathcal{L}_S \p_A = -(\p_A v^B)\p_B$)
$$\tfrac{1}{2}\varepsilon^{ABC}\,(\mathcal{L}_S\Omega)(\p_B,\p_C)
\,=\, \mathcal{L}_S \omega^A - \omega^D(\p_D v^A) + (\p_D v^D) \omega^A$$
}

\section{Volume-preserving diffeomorphisms of $\R^3$}\label{app:volpre}
\mostimportant{%
Consider the problem of constructing pairs of
(local) coordinates systems $x^A$
and $\lambda^A$ on a 3-dimensional manifold such that
$\lambda \mapsto x$ and
$x \mapsto \lambda$ are volume-preserving,
$$\dd x^1 \wedge \dd x^2 \wedge \dd x^3 
= \dd \lambda^1 \wedge \dd \lambda^2 \wedge \dd \lambda^3$$
This equation is conveniently solved
with respect to $\lambda^1,\lambda^2,x^3$ as
independent vars,
$$
\frac{\p \lambda^3}{\p x^3} = 
\frac{\p x^1}{\p \lambda^1} \frac{\p x^2}{\p \lambda^2}
- 
\frac{\p x^1}{\p \lambda^2} \frac{\p x^2}{\p \lambda^1}
$$
One can specify two functions $x^1,x^2$ freely
and then determine $\lambda^3$ by integrating.
The result is in implicit form.

When applying this to 3D incompressible Euler,
a canonical choice is to let
$\lambda^1,\lambda^2$ be two Lagrangian coordinates that
are constant along vortex lines,
a property that is preserved by the flow:
$\Omega = f(\lambda^1,\lambda^2)\,\dd \lambda^1 \wedge \dd \lambda^2$
 independent of time.
}

%%%%%%%%%%%%%%%%%%%%%%%%%%%%%%%%%%%%%%%%%%%%%%%%%%%%%%%%%%%%%%%%%%%%%
\section{Sobolev spaces on $\R^3$}\label{app:norms}
\mostimportant{We give examples for spaces $\SN$, $\SDN$
that satisfy the assumptions in Section \ref{sec:basicargXXXX}.
They depend on two arbitrary integer parameters
$\text{M} \geq 2$ and $\text{L} \geq 0$.
They are $L^2$ Sobolev spaces,
possibly with weights. The weights are such that
higher derivatives decay faster at infinity.
The spaces are subsets of the space $\mathcal{S}'$ of real-valued tempered
distributions on $\R^3$.
We are sloppy about the distinction between real-valued
and complex-valued distributions,
and use the Fourier transform $\widehat{\phantom{u}}: \mathcal{S}' \to \mathcal{S}'$.
}

{\bf Definition of $\SN$.}
Suppose $\SDN$ is a Banach space with continuous
inclusion $\SDN \subset L^2$.
Define $\SN$ in terms of $\SDN$ by
$$
\NN{u} = \begin{cases}
\infty & \text{if $\widehat{u}\notin L^1_{\text{loc}}$}\\
\NDN{\p u} & \text{if $\widehat{u} \in L^1_{\text{loc}}$}
\end{cases}$$

\mostimportant{%
This definition implies that constant functions are not in $\SN$,
except the zero function. The definition of $\SN$
is motivated by the equivalence of norms
$$
\NN{u}\;\; \sim \;\; {\textstyle\sup_B} \big(\|\widehat{u}\|_{L^1(B)} \big/ \||k|^{-1}\|_{L^2(B)}\big) + \NDN{\p u}
$$
where $B$ runs over the bounded subsets of $\R^3$.
 To see the equivalence,
if $\widehat{u} \notin L^1_{\text{loc}}$ then both sides are infinite,
and if $\widehat{u} \in L^1_{\text{loc}}$ then the
first term on the right is controlled by the second, because by
Schwarz's
inequality it is bounded by $\||k|\widehat{u}\|_{L^2} \lesssim \|\p u\|_{L^2}
\lesssim \NDN{\p u}$.
The discussion makes sense only because
the denominator is finite,
$|k|^{-1} \in L^2_{\text{loc}}$.
The space $\SN$ is a Banach space, because the right hand side of the equivalence is.
Essentially, $\NN{u}$ interpolates in Fourier space
between an $L^1$ for small $|k|$
and an $H^1$ for large $|k|$, at least.
}

\FloatBarrier
{\bf Definition of $\SDN$.}
$\NDN{u}^2 = \sum_{\alpha,\beta} \|x^{\beta}\p^{\alpha} u\|_{L^2}^2$
with summation over
$|\alpha| \leq \text{M}$
and $|\beta| \leq \text{L}$ and
$|\beta| \leq |\alpha| + 1$, see the figure for examples.

\begin{figure}[h!]
\centering
\input{norms.pstex_t}
\caption*{
\footnotesize
A box is for the sum of all $\|x^{\beta}\p^{\alpha}u\|_{L^2}$
with given $|\alpha| \geq 0$ (horizontal) and $|\beta|\geq 0$ (vertical).\\
Adding the gray box (plain $L^2$) yields an equivalent norm.
}
\end{figure}

{\bf Continuous embedding $\SN,\SDN \subset C^0$.}
Use $\|u\|_{C^0} \lesssim \|\p u\|_{L^2} + \|\p^2 u\|_{L^2}$,
a Sobolev inequality for compactly supported $u$.

{\bf Continuous embedding $\SN \subset L^p$.}
If $\text{L} = 0$,
continuous embedding iff $p \geq 6$,
by the Gagliardo-Nirenberg-Sobolev inequality.
If $\text{L} \geq 1$,
continuous embedding iff $p \geq 2$,
by $\|u\|_{L^2} \lesssim \|x\p u\|_{L^2}$
for compactly supported $u$
by partial integration.

{\bf Product inequality.}
$\NDN{uv} \lesssim \NDN{u}\NDN{v}$
by adapting the well-known
$\|fg\|_{H^m} \lesssim \|f\|_{H^2} \|g\|_{H^m}
+ \|f\|_{H^m} \|g\|_{H^2}$, $m \geq 2$ for standard $L^2$ Sobolev spaces.

{\bf Inverse of $L_0: \SN \to \SDN$.}
Let $\NN{u} < \infty$
and $\NDN{w} < \infty$.
Here $w = (w_1,w_2)$, with
vector valued
$u,w_1$ and scalar valued $w_2$,
and the integrability condition $ik \cdot \widehat{w}_1 = 0$
is understood.
Now
$$
L_0u = w \quad \Longleftrightarrow \quad 
\{k\} \widehat{u} = \widehat{w}
\quad \Longleftrightarrow \quad
\widehat{u} = \{k\}^{-1} \widehat{w}
$$
where
$\{k\}\widehat{u} = (ik\times\widehat{u}, ik\cdot\widehat{u})$
and
$\{k\}^{-1}\widehat{w} = i|k|^{-2}(k \times \widehat{w}_1-k \widehat{w}_2)$,
and below
$\{k\}^Q$ is any matrix with entries homogeneous of
degree $Q$ in $k$.
To justify the second equivalence, note that
$\widehat{u}$ and $\widehat{w}$ are in $L^1_{\text{loc}}$,
and are therefore determined by their `values' at $k \neq 0$.
It now suffices to show that
 $w \mapsto u$ with $\widehat{u} = \{k\}^{-1}\widehat{w}$
is a map $\SDN \to \SN$; it is then the unique left- and right-inverse of $L_0$.
The map gives $\widehat{u} \in |k|^{-1}L^2
\subset (L^1_{\text{loc}} \cap \mathcal{S}')$,
and one can take the inverse Fourier transform to get $u \in \mathcal{S}'$.
To see why $\NN{u} \lesssim \NDN{w}$,
take for example $(\text{M},\text{L})=(3,2)$.
\begin{figure}[h!]
\centering
\input{l0inv.pstex_t}
\end{figure}
There are separate homogeneous inequalities for {\bf a}, {\bf b}, {\bf c},
{\bf d}, {\bf e}, see the figure.
For {\bf a}, one must show
$$
\|\p(k^1\{k\}^{-1}\widehat{w})\|_{L^2}
+ \|\p^2( k^2 \{k\}^{-1}\widehat{w})\|_{L^2}
\;\lesssim\;
\|\p \widehat{w}\|_{L^2}
+ \|\p^2(k \widehat{w})\|_{L^2}
$$
with summations over all multiindices of the indicated degrees understood.
On the right hand side, we get an equivalent seminorm if
we put all the derivatives on $\widehat{w}$ directly,
but beware that this equivalence is not term by term.
The derivatives on the left produce among others
$\|\{k\}^{-1}\widehat{w}\|_{L^2}$,
which is controlled by $\|\p \widehat{w}\|_{L^2}$ by Pitt's inequality
$\||x|^{-1}f\|_{L^2} \lesssim \|\p f\|_{L^2}$.
The other terms are controlled more directly.
Similar for {\bf b}, {\bf c}, {\bf d}, {\bf e}, without Pitt's
inequality.
\FloatBarrier

%%%%%%%%%%%%%%%%%%%%%%%%%%%%%%%%%%%%%%%%%%%%%%%%%%%%%%%%%%%%%%%%%
\section{Flow map}\label{app:flowmap}

The flow map $\phi(t,s,X)$ for a vector field $V(t,X)$
is defined by
$$
\phi(t,t,X) = X
\qquad
\tfrac{\dd}{\dd s} \phi(s,t,X) = V(s,\phi(s,t,X))
$$
with the derivative at constant $t$, $X$.
With $\phi(s,t) : X \mapsto \phi(s,t,X)$ we have
$$\phi(s,u)\circ \phi(u,t) = \phi(s,t)$$
because both sides satisfy the same ordinary differential equation
in $s$, and
are the same when $s = u$. It follows that
$\phi(s,t)^{-1} = \phi(t,s)$.

\mostimportant{%
Note that
for fixed $t$,
the intervals of existence of $s \mapsto \phi(s,t,X)$
depend on $X$, say $I_X \subset \R$, and it is possible
that the intersection $\bigcap_X I_X$ contains no open neighborhood of $t$.
}

\section{Infinitesimal composition estimate}\label{app:compest}

Let $M$ be a 3-dimensional manifold.
Let $y(s): M \to \R^3$ be a diffeomorphism for every $s$,
with smooth dependence on $s$.
Let $f: M \to \R$. \emph{We estimate how the norm of
$f \circ y(s)^{-1}: \R^3 \to \R$
changes as a function of $s$, infinitesimally at $s=0$.}
The discussion is for weighted $L^2$ Sobolev norms.
The notation is local.

\mostimportant{%
We do not differentiate
$f \circ y(s)^{-1}$ with respect to $s$ and then take the norm,
because one loses derivatives that way.
Instead we differentiate the norm itself.
More precisely, we differentiate pointwise on $M$ the integrands
for weighted $L^2$ Sobolev seminorms.
}

Define the $s$-dependent 3-form
$$
\omega^{\alpha\beta} = |y^{\beta} \p^{\alpha}_yf|^2\,
\dd y^1 \wedge \dd y^2 \wedge \dd y^3
$$
Here $\p_y^{\alpha}: C^{\infty}(M) \to C^{\infty}(M)$
are the partial derivatives w.r.t.~$y$ as coordinates.
We calculate
$D\omega^{\alpha\beta}$.
Here $D$ is the derivative w.r.t.~$s$ at $s=0$.
Abbreviate $x = y(0)$ and $v = Dy$. Then
\begin{align*}
D(\dd y^1 \wedge \dd y^2 \wedge \dd y^3) & = \tr(\p_xv)\, \dd x^1 \wedge 
\dd x^2 \wedge \dd x^3\\
D y^{\beta} & = \textstyle\sum_A \beta_A x^{\beta - e_A} v^A\\
D \p^{\alpha}_y & =
-\textstyle\sum_{0 \leq \gamma < \alpha}
{\alpha \choose \gamma}
\sum_A
(\gamma + e_A)!\,
(\p_x^{\alpha-\gamma}v^A)\,\p_x^{\gamma+e_A}
\end{align*}
%------------------------------------------
where $e_A \in \N^3_0$ are the unit vectors,
and $\gamma < \alpha$ means $\gamma \leq \alpha$
and $\gamma \neq \alpha$. Then
\begin{align*}
D\omega^{\alpha\beta}
& = 2\,(x^{\beta}\p_x^{\alpha}f) b_{\alpha\beta}\, \dd x^1 \wedge \dd x^2 \wedge \dd x^3\qquad \text{with}\\
b_{\alpha\beta}
& =
(Dy^{\beta})(\p_x^{\alpha}f)
+ x^{\beta} (D \p_y^{\alpha})f
+ \tfrac{1}{2}\,x^{\beta}(\p_x^{\alpha}f)\tr(\p_xv)
\end{align*}
\noindent
{\bf Remark 1:}
Each term in $b_{\alpha\beta}$ is, up to constant factors,
and with $\p = \p_x$,
\begin{itemize}
\itemsep 0pt
\item either bounded pointwise by $\langle x \rangle^{|\beta|} |\p^{\alpha}f|\,|\langle x \rangle^{-1} v|
$ with $\langle x \rangle = (1+|x|^2)^{1/2}$,
\item or of the form $x^{\beta}(\p^{\mu}(\p f))(\p^{\nu}(\p v))$
with $|\mu|+|\nu| = |\alpha|-1$,
\item (only for $\alpha = 0$) or of the form $x^{\beta}f\,(\p v)$.
\end{itemize}

\mostimportant{In the first item,
it may be convenient to use
$\|\langle x \rangle^{-1} v \|_{\infty}
\lesssim |v(0)| + \|\p v\|_{\infty}$.}

%------------------------------
\noindent{\bf Remark 2:} Formally by Cauchy-Schwarz,
for any finite set $S$ of pairs $(\alpha,\beta)$,
$$
\big|
D
\big(\textstyle\sum_{(\alpha,\beta) \in S}
\int_M \omega^{\alpha\beta}
\big)^{1/2}
\big|
 \leq \big(\textstyle\sum_{(\alpha,\beta) \in S}
\|b_{\alpha\beta}\|_{L^2(\dd^3x)}^2\big)^{1/2}
$$

\mostimportant{%
Different $S$ yield different weighted Sobolev seminorms.
One can use Remark 1 with a bilinear estimate
and Remark 2 to obtain estimates of the form
$|D\|f \circ y(s)^{-1}\|_A| \lesssim \|v\|_B\,\|f\|_A$.
}

%%%%%%%%%%%%%%%%%%%%%%%%%%%%%%%%%%%%%%%%%%%%%%%%%%%%%%%%%%%%%%%%%5
\section{Continuous self-similarity}\label{app:css}
\mostimportant{We write down the
equations for continuously self-similar
solutions to incompressible Euler,
in geometric form analogous to \eqref{eq:geomform}.
The role of the constant $\KA > 0$ is to select
 a 1-dimensional subgroup of the
2-dimensional scaling symmetry group of incompressible Euler.
}

On a 3-dimensional manifold $M^3$, consider the system
\begin{subequations}\label{eq:css}
\begin{align}
\mathcal{L}_Z \xi^A & = \zeta^A + \KA\, \xi^A\\
\mathcal{L}_Z \Gamma & = (2\,\KA - 1) \Gamma \label{eq:euzehbdfd}\\
\dd \Gamma & = 0\\
\delta_{AB}\, \dd \zeta^A \wedge \dd \xi^B & = \Gamma\\
\varepsilon_{ABC}\, \dd \zeta^A \wedge \dd \xi^B \wedge \dd \xi^C & = 0
\end{align}
\end{subequations}
A solution to \eqref{eq:css} yields a self-similar solution to \eqref{eq:geomform}
on $t < 0$, as follows. Let $M^4 = \{t < 0\} \times M^3$.
The objects in \eqref{eq:css}
induce objects on $M^4$ that do not depend on $t$,
and that we denote by the same symbols. Now
\begin{align*}
S & = \p_t + |t|^{-1} Z&
x^A & = |t|^{\KA} \xi^A&
v^A & = |t|^{\KA-1} \zeta^A&
\Omega & = |t|^{2\,\KA-1} \Gamma
\end{align*}
solves \eqref{eq:geomform} by direct calculation.

\mostimportant{%
We have assumed that $(\xi^1,\xi^2,\xi^3): M^3 \to \R^3$ is a diffeomorphism.
Observe that if one uses $t,x^A$ as independent variables,
then $v^A(qt,q^{\KA}x) = q^{\KA-1} v^A(t,x)$ for all $q > 0$.
}

%%%%%%%%%%%%%%%%%%%%%%%%%%%%%%%%%%%%%%%%%%%%%%%%%%%%%%%%%%%%%%%%%%
\section{Like $x^{\alpha m+ n}$ but linearly independent}\label{app:specialfunc}

\mostimportant{Consider for $\alpha > 0$ the family of functions
$\mathcal{F}_{\alpha} = \{x \mapsto x^{\alpha m + n}\}_{m,n \geq 0}$
indexed by integers $m,n$. Here $x > 0$.
The family $\mathcal{F}_{\alpha}$ is linearly independent iff
$\alpha \notin \Q$, with $\Q$ the rationals.
We apply a meromorphic in $\alpha$ (with poles in $\Q$)
change of basis to $\mathcal{F}_{\alpha}$ to obtain a new family that
does not suffer from loss of linear independence.}

\mostimportant{%
We do not use logarithms to
resolve the dependencies `piecewise in $\alpha$'.
Instead, the family constructed below is smooth in $\alpha$.
This is useful
for applications in which
$\alpha$ is an unknown,
say to differentiate with respect to $\alpha$ in a Newton iteration.

Since $x^{\alpha m+n} = y^{(1/\alpha)n+m}$
with $y = x^{\alpha}$, the construction yields two families.
}

{\bf Construction.}
Fix a compact interval $A \subset (0,\infty)$.
Set
$\mathbf{s}_k = (-1)^k$ except for $\mathbf{s}_0 = \tfrac{1}{2}$.
Define  $(g_m)_{m \geq 0}$ in terms of $(f_\ell)_{\ell \geq 1}$
by a generating function:
\begin{subequations}
\begin{align}\label{euiehrff}
\textstyle\sum_{m \geq 0} t^m\,g_m
&= \exp \big(\textstyle\sum_{\ell \geq 1} t^\ell f_\ell\big)\\
\label{aswq}
f_\ell
&=
\frac{i\alpha}{\pi}
\bigg(
\sum_{k \geq 0}  \frac{\mathbf{s}_k x^{\alpha \ell}}{(\alpha \ell)^2-k^2}
- \sum_{k \in \Z \cap (\ell A)}  \frac{\mathbf{s}_k x^k}{(\alpha \ell)^2-k^2}\bigg)
\end{align}
\end{subequations}
with $i = \sqrt{-1}$.
The $f_{\ell}, g_m$
are polynomials jointly in $x^{\alpha}$ and $x$,
with complex coefficients that
are meromorphic in $\alpha$
with poles in $\mathbbm{Q}$,
but the $f_{\ell}, g_m$ themselves are holomorphic for $\alpha \in A$
since the poles in $A \cap \mathbbm{Q}$ cancel.
The
family of complex functions
$\{g_m(x)x^n\}_{m,n \geq 0}$
is
linear independent for $\alpha \in A$.
\step

{\bf Remarks.}
$$
\begin{pmatrix}
g_0\\
g_1\\
g_2\\
\vdots
\end{pmatrix}
=
\begin{pmatrix}
\ast_0 & 0 & 0 & \ldots\\
\ast_1 & \ast_0 & 0 & \ldots\\
\ast_2 & \ast_1 & \ast_0 & \ldots\\
\vdots & \vdots & \vdots & \ddots
\end{pmatrix}
\begin{pmatrix}
x^{0\alpha}\\
x^{1\alpha}\\
x^{2\alpha}\\
\vdots
\end{pmatrix}
$$
Each $\ast_k$ is a polynomial in $x$,
with monomial degrees in $kA$,
and with coefficients meromorphic in $\alpha$.
The diagonal entries of degree 0
\emph{do not vanish for any $\alpha \in \C$}:
\begin{equation*}
\text{(coefficient of $x^{m \alpha}$ in $g_m$)} =
2^{-m} \mathbf{c}_m \big/ \textstyle\prod_{j=1}^m \sin (\pi \alpha j)
\end{equation*}
where $\mathbf{c}_m = i^m \exp(-\pi i \alpha m(m-1)/2)$ has unit length
for real $\alpha$.

\mostimportant{
To check this,
one can use the following
partial fraction and Pochhammer identities:
$$
\frac{1}{2z}\, \frac{\pi}{\sin(\pi z)}
= \sum_{k \geq 0} \frac{\mathbf{s}_k}{z^2-k^2}
\qquad\qquad
\sum_{m \geq 0} \frac{t^m}{\prod_{j=1}^m (1-z^j)} = \exp\Big(\sum_{\ell \geq 1}\frac{1}{\ell}\,\frac{t^\ell}{1-z^\ell}\Big)
$$}

Set $h_{\ell} = (x\p_x - \alpha \ell) f_{\ell}$.
Applying $x \p_x - \alpha t \p_t$ to
\eqref{euiehrff}
gives
$$(x\p_x - \alpha m) g_m = \textstyle \sum_{\ell = 1}^m h_{\ell} g_{m-\ell}$$
This is a recursion for the $g_m$ given the $h_{\ell}$,
up to constants of integration.
The $h_{\ell}$ are polynomials in $x$, with coefficients
that are holomorphic in $\alpha \in A$.

\mostimportant{Explicitly
$h_{\ell} = (i\alpha/\pi) \sum_{k \in \Z \cap (\ell A)}
(-1)^k x^k/(\alpha \ell + k).$}

We have $g_m g_n = \sum_{0 \leq k \leq m+n} p_{mnk} g_k$
and
$(x\p_x) g_m = \sum_{0 \leq k \leq m} q_{mk} g_k$
where $p_{mnk}$ and $q_{mk}$ are polynomials in $x$,
with monomial degrees in
$(m+n-k)A$ and $(m-k)A$ respectively,
and with coefficients holomorphic in $\alpha \in A$.
They are `differential algebra' structure coefficients.

The family of real functions
$\{\text{Re} (g_m(x)/\mathbf{c}_m) x^n\}_{m,n \geq 0}$
has the same linear span, and similar properties.

\mostimportant{Note that for real $\alpha$, the $f_{\ell}$ are purely
imaginary, and $(\sum_m t^m g_m)(\sum_m t^m \overline{g_m}) = 1$.}

%%%%%%%%%%%%%%%%%%%%%%%%%%%%%%%%%%%%%%%%%%%%%%%%%%%%%%%%%%%%%%%%%%
\section{Stabilizer candidates}\label{app:somestabcand}

\mostimportant{
Suppose a group $G$ acts on the
set of solutions to incompressible Euler
(symmetry of the equations).
Let $H \subset G$ be the stabilizer of some solution
(symmetry of a solution).
The candidates for $H$ can be obtained by
classifying the subgroups of $G$, up to conjugacy.

One can take $G$ to be the Galilean group and two scalings,
see Appendix \ref{app:sym}. We are much less thorough
and consider only a smaller group,
and only a partial classification.
}

Set
$G = G_1\times G_2 \times G_3$ with
$G_1 = \R^+ \times \R^+$
and
$G_2 = \{\pm 1\}$ 
and
$G_3 = \SO(3)$.
Then $G_2 \times G_3 = \OO(3)$.
The group $G$ acts faithfully on velocities $v(t,x)$ by
\begin{align*}
\big(((a,b),\sigma,R)v\big)(t,x)
& = (a/b)\, \sigma R\, v(at, b \sigma R^{-1} x)
\end{align*}
%-----------------------------------------------------------
Let $H \subset G$ be a subgroup.
Let $H_i = p_i(H)$ with
$p_i: G \to G_i$ the projections.
We only consider $H$ for which the $H_i$ are, up to conjugacy in $G_i$:
\begin{alignat*}{6}
H_1 &\qquad & \R^+ & \to G_1\qquad & a &\mapsto (a,a^{\alpha}) \\
    & & \Z   & \to G_1 & n &\mapsto (a_1^n,(a_1^n)^{\alpha}) \\
\rule{0pt}{13pt} H_2 && \{1\} & \to G_2 & 1 & \mapsto 1\\
    && \{\pm 1\} & \to G_2 & \sigma & \mapsto \sigma\\
\rule{0pt}{13pt} H_3 && \{1\} & \to G_3 & 1 & \mapsto \mathbbm{1}_{3 \times 3}\\
    && \SO(2) & \to G_3 & r & \mapsto \mathbbm{1}_{1 \times 1} \oplus r\\
    && \OO(2) & \to G_3 & r & \mapsto (\det r) \oplus r
\end{alignat*}
%--------------------------------------------------------------
\mostimportant{%
Changing the parameters $\alpha > 0$, $a_1 > 1$ yields
non-conjugate subgroups of $G_1$, because it is Abelian.
We use $\oplus$ to build block diagonal matrices.
We have selected subgroups as follows:
For $H_1\subset G_1$ only closed `1D subgroups'.
For $H_3\subset G_3$ we have excluded
discrete (i.e.~finite) nontrivial subgroups,
$\SO(3)$ itself,
but also
non-closed subgroups such as the infinite subgroups that are generated
by a single element of $\SO(3)$.
Since $H_3$ is the image of a projection,
we have no good justification for taking only closed $H_3$.
%---------------------------------------------------------------

Note that
$H \subset H_1 \times H_2 \times H_3$
is a subdirect product.
Set $H_{23} = p_{23}(H)$ with $p_{23}: G \to G_2\times G_3$.
Then $H \subset H_1 \times H_{23}$
and $H_{23} \subset H_2 \times H_3$ are also subdirect products. 
A subdirect product of two groups is a fiber product, by Goursat's lemma.
We first determine all $H_{23}$, then all $H$.
For $H_{23}$, consider the 6 pairs $(H_2,H_3)$ separately.
There are 7 cases for $H_{23}$,
\begin{itemize}
\item the 6 trivial products $H_{23} = H_2 \times H_3$,
\item with $H_2 = \{\pm 1\}$, $H_3 = \OO(2)$
the nontrivial $H_{23} = \{(\sigma,r) |\,\sigma = \det r\} \cong \OO(2)$.
\end{itemize}
For $H$, there are 14 pairs $(H_1,H_{23})$ to consider.
The results are below.}

The group $H \subset H_1 \times H_2 \times H_3$ is one of the following:
\begin{center}
\begin{tabular}{l|l|l|l|l|l|l}
& $H_1$ & $H_2$ & $H_3$ & equations that define $H$ & dim & no swirl\\
\hline
1 & any & any & any & none & 0 or 1 or 2 & y or n\\
2 & any & $\{\pm 1\}$ & $\OO(2)$ & $\sigma = \det r$ & 1 or 2 & y\\
3 & $\R^+$ & any & $\SO(2)$ & $\mu \log a = \ANGLE(r)$, $\mu\neq 0$ & 1 & n\\
4 & $\Z$ & $\{\pm 1\}$ & any & $(-1)^n = \sigma$ & 0 or 1 & n\\
5 & $\Z$ & any & $\OO(2)$ & $(-1)^n = \det r$ & 1 & n\\
6 & $\Z$ & $\{\pm 1\}$ & $\OO(2)$ & $(-1)^n = \sigma (\det r)$ & 1 & y\\
7 & $\Z$ & $\{\pm 1\}$ & $\OO(2)$ & $(-1)^n = \sigma = \det r$ & 1 & n
\end{tabular}
\end{center}
\mostimportant{Row 1 are the 12 direct products
$H = H_1 \times H_2 \times H_3$.
`No swirl' means that $H$ contains the $\OO(2)$
subgroup of $\{\text{identity}\}\times G_2 \times G_3$
given by $(\sigma,r) \in \{\pm 1\}\times \OO(2)$ with
$\sigma = \det r$.}
%%%%%%%%%%%%%%%%%%%%%%%%%%%%%%%%%%%%%%%%%%%%%%%%%%%%%%%%%%%%%%%%%%%%%%%
\end{document}